\documentclass{amsart}

\newtheorem{propo}{Proposition}[section]

\newtheorem{lemma}[propo]{Lemma}
\newtheorem{corol}[propo]{Corollary}

\newtheorem{theo}[propo]{Theorem}

\newtheorem{remar}[propo]{Remark}

\newcommand{\bl}{\begin{lemma}\label}
\newcommand{\el}{\end{lemma}}

\newcommand{\ld}{,\ldots ,}

\newcommand{\ra}{ \rightarrow }

\newcommand{\se}{ \subseteq }
\newcommand{\lan}{ \langle }
\newcommand{\ran}{ \rangle }

\newcommand{\diag}{\mathop{\rm diag}\nolimits}

\newcommand{\Id}{\mathop{\rm Id}\nolimits}

\newcommand{\FF}{\mathop{\mathbb F}\nolimits}

\newcommand{\ZZ}{\mathop{\mathbb Z}\nolimits}

\newcommand{\al}{\alpha}
\newcommand{\be}{\beta}

\newcommand{\ep}{\varepsilon}
\newcommand{\lam}{\lambda }

\newcommand{\om}{\omega }
\newcommand{\Om}{\Omega }

\newcommand{\up}{^{-1}}

\newcommand{\si}{\sigma }
\newcommand{\med}{\medskip}
\newcommand{\bp}{\begin{proof}}
\newcommand{\enp}{\end{proof}}

\def\d12{{_{12}}}

\def\acf{{algebraically closed field }}

\def\ei{{eigenvalue }}
\def\eis{{eigenvalues }}

\def\f{{following }}

\def\ho{{homomorphism }}

\def\ii{{if and only if }}

\def\ir{{irreducible }}

\def\irr{{irreducible representation }}

\def\st{Suppose that }
\def\hw{highest weight }

\def\itf{{It follows that }}

\def\rep{{representation }}
\def\reps{{representations }}

\def\st{{Suppose that }}

\begin{document}

\title[rational elements of algebraic groups]{Rational elements in representations of simple algebraic groups, I}

\author[A. Zalesski]{Alexandre Zalesski}
\address{School of Mathematics, University of East Anglia, Norwich, UK }

\thanks{e-mail:
   alexandre.zalesski@gmail.com}

\maketitle

\centerline {\it Dedicated to Pham Huu Tiep on the occasion of his $60$-th birthday}

\def\ag{algebraic group }

\bigskip
{\small{\bf Abstract}   An element $g$ of a group $G$ is called {\it rational} if $g$ is conjugate to $g^i$ for every integer $i$ coprime to  $|g|$. We determine all triples $(G,g,\phi)$, where $G$ is a simple algebraic group of type 
$A_n,B_n$ or $C_n$ over an algebraically closed field of characteristic $p\geq 0$, $g\in G$ is a rational  odd order semisimple element and $\phi$ is an irreducible representation of $G$  such that $\phi(g)$ has eigenvalue 1.
}



\section{Introduction}

The primary motivation of this work lies in the study of \eis of semisimple elements of simple algebraic groups 
in their \ir representations.  This problem is not treatable in full generality,
however, a number of special cases of it  were discussed in literature. The earliest work is probably by Wilson \cite{Wi} who determined the \ir \reps of simple algebraic groups in characteristic 2 in which some element of order 3 does not have \ei 1.  In \cite{Z91} there  appeared a number of useful observations, and certain comments are available in surveys \cite{Z88} and \cite{z09}. The occurrence of the \ei 1 deserves a particular attention, as this carries some geometrical meaning, and is significant for applications. Robert Guralnick and Pham Huu Tiep \cite{GT} determine the simple algebraic groups $G$ whose every element in {\it every}  \rep has \ei 1. (If $\rho$ is a \rep of $G$ then every element of $G$ has weight 0 \ii $\rho$ has weight 0; the simple  simply connected groups whose all 
\reps have weight 0 were listed in \cite{Z91} and the groups of adjoint type with this property were determined in  \cite{GT}. In addition, in \cite{GT} the authors determined all finite groups of Lie type with the above \ei property, some partial results appeared in \cite{Z91}.)   In \cite{Z17,Z20} 
 for groups $G=SL_n(2)$ and $Sp_{2n}(2)$ I determine the \ir 2-modular \reps of $G$ in which every element of $G$ has \ei 1.  It is clear that one cannot expect to determine the elements $g\in G$ which have \ei 1 in every given \ir representation. 
In this paper  we consider  elements of "highest symmetry",
specifically, the rational  elements.  Recall that an element $g$ of a group $G$ is called {\it rational} if $g$ is conjugate to $g^i$ for every integer $i$ coprime to  $|g|$. In other words, if $\lan g\ran$ is  the group generated by $g$ then all generators of $\lan g\ran$ are conjugate in $G$, equivalently,  $N_G(\lan g\ran)/C_G(\lan g\ran)\cong {\rm Aut}\,(\lan g\ran)$. This implies that for every \rep $\phi$  of $G$ and every \ei $\lam$ of $\phi(g)$ the powers $\lam ^i$  are \eis of $\phi(g)$ for all $i$ coprime to $|g|$.

Let $G=G(F)$ be a simple simply connected  \ag over an \acf $F$ of characteristic $p\geq 0$.  
In Theorems \ref{th1}, \ref{th3} and \ref{th2} below   we determine all triples $(G,g,\phi)$, where $G$ is
of type $A_n,B_n$, $C_n$, respectively, in 
$g\in G$ is a rational  odd order semisimple element and $\phi$ is an \irr of $G$  such that $\phi(g)$ has \ei 1. Simple algebraic groups of the other   types will be considered in a subsequent paper. With additional efforts, one can probably drop the restriction that $|g|$ is odd. Note that $\lan g \ran\cap Z(G)=1$ whenever $g$ is rational and $|g|$ is odd.

\med

To state our results,
 we use a parametrization of the \ir \reps of $G$ in terms of highest weights. If the rank of $G$ equals $n$ then we denote by $\om_1\ld \om_n$ the fundamental weights (ordered as in \cite{Bo}).
If $p=0$ then $p^t\om_i$ is interpreted as $\om_i$. To describe the elements $g\in G$ in question, we use their Jordan form  at the natural \rep of $G$, which we define as the \irr of $G$ with highest weight $\om_1$. An \irr of $G$ with highest weight $\om$ is denoted by $\rho_{\om}$ unless otherwise stated. 

Given a natural number $m$, we denote by $R(m)$ and $\Phi(m)$ the set of all $m$-roots of unity and all primitive $m$-roots of unity, respectively, and set $\phi(m)=|\Phi(m)|$. The mapping $m\ra \phi(m)$ for $m\in \ZZ^+$ is called the Euler function.
For a \rep $\rho$ of $G$ denote by $E(\rho(g))$ the set of all \eis of $\rho(g)$ disregarding the multiplicities. If $g$
is of order $m$ then $E(\rho(g))\subseteq R(m)$.

\begin{theo}\label{th1} Let $G=SL_n(F), n>1,$ and let $g\in G$ be a semisimple odd order  element. Let $\rho$ be an \irr of G with \hw $\om$. Suppose that $g$ is rational. Then $1$ is an \ei of $\rho(g)$ unless one of the \f holds for some integer $t\geq 0$:

\med
$(1)$ $\dim\rho=n$, n even, and $\om\in\{p^t\om_1,p^t\om_{n-1}\};$

$(2)$ $n=6$, $|g|=9$ and $\om=p^t\om_3;$

$(3)$ $n=6$, $|g|=15$, $g=\diag(g_1,g_2)$, $|g_1|=5,|g_2|=3$ and $\om=p^t\om_3;$

$(4)$ $n=8$, $|g|=15$ and $\om\in\{p^t\om_3,p^t\om_5\};$

$(5)$ $n=8$, $|g|=15$, $g=\diag(g_1,g_2,g_3)$, $|g_1|=5,|g_2|=|g_3|=3$ and  $\om\in\{p^t\om_3,p^t\om_5\};$

$(6)$ $n=10$, $|g|=45$,  $g=\diag(g_1,g_2)$, $|g_1|=5,|g_2|=9$ and
 $\om\in\{p^t\om_3,p^t\om_7,p^t\om_5\};$

$(7)$ $n=14$, $|g|=45$,  $g=\diag(g_1,g_2)$, $|g_1|=15,|g_2|=9$ and   $\om\in\{p^t\om_3,p^t\om_{11}\}.$
\end{theo}

Here $g_i$, $i=1,2,3$ is a diagonal matrix of size $\phi(|g_i|)$ with diagonal entries $\Phi(|g_i|)$, and $g$ is identified with $\rho_{\om_1}(g)$. 

As a by-product of our proof of Theorem \ref{th1}, we have

\begin{theo}\label{tt9} In assumption of  Theorem  {\rm \ref{th1}}, every \ei of $\rho_{\om_1}(g)$ is an \ei of $\rho(g)$, unless possibly  for the case where $n$ is even, $\om=\om_i$ with $i$ even, $|g|=3l$ with $(3,l)=1$.
\end{theo}

A more precise description of the exceptional case in  Theorem  \ref{tt9}  is given in Lemma \ref{un2}.

\begin{theo}\label{th3} Let $G=Spin_{2n+1}(F), n>1$, $p\neq 2$, and let $g\in G$ be a semisimple  element. Let $\rho$
be an \irr of G with \hw $\om$. Suppose that $g$ is rational in G and $|g|$ is odd. Then $1$ is an \ei of $\rho(g)$ unless one of the \f holds for some integer $t\geq 0$:

\med
$(A)$ $\om=p^t\om_n$  and  one of the \f holds:

\med
$\,\,\,\,\,\,\,\,\,\,~~~~~~~~~~~~~~(1)$ $E(\rho(g))=\Phi(m)R(|g|/m)$ with $m\in\{3,5,9\}$ and $(m,|g|/m)=1;$ in other words, $E(\rho(g))$\\ \hskip1cm consists of all  $|g|$-roots of unity whose orders are multiples of $m;$

$\,\,\,\,\,\,\,\,\,\,~~~~~~~~~~~~~~(2)$ 
$E(\rho(g))=\Phi(5)(R(m)\setminus 1) R(|g|/5m)$ with $m\in\{3,9\}$ and $(5m,|g|/5m)=1;$
in other words, $E(\rho(g))$ consists of all $|g|$-roots of unity whose orders are multiples of $15;$

\med
$(B)$ $\om=p^a\om_1+p^b\om_n$, $a,b\in \ZZ^+ $, and g is as in $(A)(2)$ above.
\end{theo}

\med
The Jordan normal form of $g$ on the natural $FG$-module can chosen to be 
 $D=\diag(D_1\ld D_r)$, where $D_i$, $i=1\ld r,$ is a diagonal matrix of size $\phi(m_i)$ with non-zero entries $\Phi(m_i)$  (for some integer $m_i\geq 1$). 
 Then $(2)$ holds \ii there is exactly one $m_i\in\{3,5,9\}$   and $(m_i,m_j)=1$ for  $i,j\in\{1\ld r\}  $, $j\neq i$, and  $(3)$ holds \ii there is exactly one $m_i\in\{3,9\}$,   exactly one $j$ with $m_j=5$  and $(15,m_k)=1$ 
for $i,j,k\in\{1\ld r\}$,  for $k\neq i , j$.

 \begin{theo}\label{th2}  Let $G=Sp_{2n}(F)$,  and let $\rho$ be an \irr of G with \hw $\om$. Let $g\in G$ be a rational element of odd order. Suppose that $\rho(g)$ does not have \ei $1$. Then one of the \f holds:

\medskip
$(A)$ $\om=p^t\om_1;$

\medskip
$(B)$ $\om=p^t\om_3$ and one of the \f holds:

\medskip
$\,\,\,\,\,\,\,\,\,\,~~~~~~~~~~~~~~(1)$ $n=3$, $|g|=9;$ 

$\,\,\,\,\,\,\,\,\,\,~~~~~~~~~~~~~~(2)$ $n=3$, $|g|=15$, $g=\diag(g_1,g_2)$, $|g_1|=5,|g_2|=3;$ 

$\,\,\,\,\,\,\,\,\,\,~~~~~~~~~~~~~~(3)$ $n=4$, $|g|=15;$ 

$\,\,\,\,\,\,\,\,\,\,~~~~~~~~~~~~~~(4)$ $n=4$, $|g|=15$, $g=\diag(g_1,g_2,g_3)$, $|g_1|=5,|g_2|=|g_3|=3;$ 

$\,\,\,\,\,\,\,\,\,\,~~~~~~~~~~~~~~(5)$ $n=5$, $|g|=45$,  $g=\diag(g_1,g_2)$, $|g_1|=5,|g_2|=9;$ 

$\,\,\,\,\,\,\,\,\,\,~~~~~~~~~~~~~~(6)$ $n=7$, $|g|=45$,  $g=\diag(g_1,g_2)$, $|g_1|=15,|g_2|=9;$

\medskip
$(C)$ $n = 5$, 
$\om=p^t\om_5$, $g=\diag(g_1,g_2)$, $|g_1|=5$, $|g_2|=9;$   

\medskip
$(D)$  $p=2$, $\om=2^t\om_n$ and one of the \f holds: 

\medskip
$\,\,\,\,\,\,\,\,\,\,~~~~~~~~~~~~~~(1)$ $g=\diag(g_1,y)$, where $g_1\in Sp_{2k}(F)$, $k=1,2,3$, $|g_1|=2^k+1$ and
 $y\in Sp_{2n-2k}(F)$ is an arbitrary rational element with $(|y|,|g_1|)=1;$

$\,\,\,\,\,\,\,\,\,\,~~~~~~~~~~~~~~(2)$ $g=\diag(g_1,g_2,y)$, where $g_i\in Sp_{2k_i}(F)$, $(k_1,k_2)=(1,2)$ or $(2,3)$,
$|g_i|=2^{k_i}+1$ for $i=1,2$  and $y\in Sp_{2n-2(k_1+k_2)}(F)$ is an arbitrary rational
element with $(|y|,|g_1g_2|)=1;$

\medskip
$(E)$   $p=2$,  $\om=2^s\om_1+2^t\om_n$, $s,t\geq 0$,  $g=\diag(g_1,g_2,g_3)$, and one of the \f holds:

\medskip
$\,\,\,\,\,\,\,\,\,\,~~~~~~~~~~~~~~(1)$  
$g_1\in Sp_{4}(F)$, $|g_1|=5$, $g_2\in Sp_{2}(F)$, $|g_2|=3$, $g_3\in Sp_{2(n-3)}(F)$  is an arbitrary rational element of order coprime to $15;$

$\,\,\,\,\,\,\,\,\,\,~~~~~~~~~~~~~~(2)$ $g_1\in Sp_{4}(F)$, $|g_1|=5$, $g_2\in Sp_6(F)$, $|g_2|=9$ and
 $g_3\in Sp_{2(n-5)}(F)$ is an arbitrary rational element
of order coprime to $15$. 
 \end{theo}

If $p\neq 2$ and $\dim \rho>1$ then $E(\rho_{\om_1}(g))\subset E(\rho(g))$ for $G=Sp_{2n}(F)$; see Lemma \ref{od1}. As above, $g$ is identified here with $\rho_{\om_1}(g)$

Other results of  the paper refine the above results for some special cases. In particular, in some cases we observe that   $E(\rho(g))=R(|g|)$. Note that with some further analysis one can hope to determine all cases where $E(\rho(g))=R(|g|)$.

 \med

{\it Notation}\,\, 
$\ZZ$ is the ring of integers, and $\ZZ^+$ is the set of non-negative integers. For $a,b\in \ZZ$, $ab\neq 0$, we write $(a,b)$ for the greatest common divisor of $a,b$, and $a|b$ means that $b$ is a multiple of $a$. A diagonal $(n\times n)$-matrix with sequential entries $x_1\ld x_n$ is denoted by $\diag(x_1\ld x_n)$.  By $\FF_q$ we denote the finite field of $q$ elements. 
All \reps below are over an \acf $F$ of characteristic $p\geq 0$ (unless otherwise   stated).
For a semigroup $H$ and subsets $X,Y\subset H$ we write $XY$ or $X\cdot Y$ for the set $\{xy:x\in X, y\in Y\}$.

For $m\in \ZZ^+$ we denote by $R(m)$ the set of all $m$-roots of unity, and by $\Phi(m)$ the set of all primitive $m$-roots of unity. So $|R(m)|=m$ and $|\Phi(m)|=\phi(m)$; we call $\phi$  the Euler function (some authors use the term "the Euler 'totient' function").
Notations  $\Lambda_i(\Phi(m))$ and $\Delta(\Phi(m))$ are introduced in Section 2.

If $\phi$ is a \rep of a group $X$ and $Y\subset X$ is a subgroup then $\phi|_Y$ means the restriction of $\phi$ to $Y$.

Notation for algebraic groups and their weight systems will be introduced in Section 2.



\def\se{semisimple }

    \section{Preliminaries}

\subsection{Cyclic groups and the  Euler function}

In this section we discuss some properties of the set $\Phi(m)$ with $m\in \ZZ^+$ which are useful 
in what follows.  

\bl{t11} Let $m$ be an odd natural number. Then $\Phi(m)\cdot \Phi(m)=R(m)$. \el

\bp   Suppose first that  $m$ is a $p$-power  for some prime $p$. Then $ R(m)=R(m')\cup \Phi(m)$, where $m'=m/p$. Note that $ R(m)$ is a group and $ R(m')$ is a subgroup of index $p$. Let $\eta\in \Phi(m)$. Then $\Phi(m)=\cup_{i=1}^{p-1} \eta^i R(m')$ and  $\eta\cdot \eta^iR(m')=\eta^{i+1}R(m')$, whence $R(m')$ and $\eta^jR(m')\subset \Phi(m)\cdot \Phi(m)$ for $j\neq 1$.
Also $\eta R(m')=\eta^2\cdot \eta^{p-1}R(m')$, whence the result.

Let   $m=m_1\cdots m_k$, where $k>1$ and $m_1\ld m_k$ are prime powers coprime to each other.
Then $\Phi(m)=\Phi(m_1)\cdots \Phi(m_k)$ and hence $\Phi(m)\cdot \Phi(m)=\Phi(m_1)\cdot \Phi(m_1)\cdots \Phi(m_k)\cdot \Phi(m_k)=R(m_1)\cdots R(m_k)=R(m)$. \enp

\bl{t2p} Let $m_1,m_2$ be odd natural numbers   and $l=(m_1,m_2)>1$. Let  $M_i=R(m_i)\setminus 1$
for $i=1,2$.  Let $m_3=m_1m_2/l$ be the least common multiple of $m_1,m_2$.
 Then $ M_1M_2=R(m_3)$.\el

\bp By Lemma \ref{t11}, the lemma holds if $m_1=m_2$, and, in addition, $R(l)\subseteq M_1M_2$.
Clearly,  $ M_1M_2$ contains 1 and $R(m_1)R(m_2)\setminus (R(m_1)\cup R(m_2))$. So it suffices to show that $R(m_1), R(m_2)\subseteq M_1M_2$. By assumption, $\Phi(l)\subseteq M_1$, so there are   distinct $\zeta_1,\zeta_2\in \Phi(l)$. Then $\zeta_1M_2\cup \zeta_2M_2=R(m_2)\subseteq M_1M_2$. Similarly, $R(m_1)\subseteq M_1M_2.$\enp

For $i\leq \phi(m)$ set $\Lambda_i(\Phi(m))=\eta_1\cdots \eta_i$, where $\eta_1,\ldots ,\eta_i\in\Phi(m)$ are distinct.
In other words,  $\Lambda_i(\Phi(m))$   is the set of all elements which are products of $i$ distinct elements in $\Phi(m)$.

\begin{lemma}\label{b33}
 For $\zeta\in R(m)$, m odd,  there exists $\eta\in\Phi(m)$ such that $\zeta\eta\in\Phi(m)$, and $\eta\neq\zeta$ if $m>3$. \end{lemma}

\bp Note that $\Phi(m)=\Phi(m_1)\cdots\Phi(m_k)$, where $m_1\ld m_k$ are prime powers  coprime to each other. Suppose first that $k=1$. If $m>3$ is a prime then the claim is obvious. Otherwise, let $m=m'p$ for $p$ a prime. 
If $\zeta\in R(m')$ then the claim is true for any $\eta\in \Phi(m)$. Otherwise, $\zeta\in\Phi(m)$.  If $p>3$ then take $\eta=\zeta^2$. If $p=3$ then we can take $\eta=\zeta^4$  if $m>3$ and $\eta=\zeta$ for $m=3$. 

Let $k>1$. If $3\notin\{m_1\ld m_k\}$ then the result follows from the above. Otherwise we can assume that $m_1=3$.
Let $m'=m/3$. Then either $\zeta\in R(m')$ and we take $\eta=\eta_3\eta'$, where $\eta_3\in\Phi(3)$, $\eta'\in\Phi(m')$ is such that $\zeta \eta'\in\Phi(m')$, or $\zeta=\eta^{\pm 1}_3\zeta'$ with $\zeta'\in R(m')$ and 
then we take $\eta=\eta^{\pm 2}_3\eta'$, where $\eta'\in\Phi(m')$ is such that $\zeta' \eta'\in\Phi(m')$, which is possible by the above. \enp

\bl{33a} Let $ m> 3$ be an odd natural number.
Then $\Lambda_2(\Phi(m))=R(m)$. In other words,
for every $\zeta\in R(m)$ there are distinct roots $\eta_1,\eta_2\in \Phi(m)$
such that $\eta_1\eta_2=\zeta$. 
\end{lemma}

\bp 
By Lemma \ref{b33}, $\zeta=\eta_1\up\eta_2$ for some $\eta_1,\eta_2\in\Phi(m)$, and $\eta_1\up\neq \eta_2$ whenever $|\zeta|<m$. If  $|\zeta|=m$ then choose $\eta_{1}=\zeta^2,\eta_2=\zeta\up$. \enp

\bl{333} Let $ m\neq 1,3,5,9,15$ be an odd natural number.
 Then there are six distinct roots  $\eta_{1}^{\pm 1},\eta_{2}^{\pm 1}, \eta_{3}^{\pm 1}\in \Phi(m)$ such that $\eta_{1}\eta_2\eta_{3}=1$.  Equivalently, for every $\eta\in \Phi(m)$ there exist $\eta_1,\eta_2\in\Phi(m)$
 such that $\eta=\eta_1\eta_2$ and $\eta^{\pm1},\eta_{1}^{\pm 1}, \eta_{2}^{\pm 1}$ are distinct.
\el

\bp
 Let $\eta\in\Phi(m)$. Note that $\eta^{2^i}\in\Phi$ for every $i\in\ZZ^+$. 

\med
(i) The lemma is true if $(m,3)=1$ and $m>5$.

\med
Indeed,
$\eta\cdot\eta^2=\eta^3\in\Phi(m)$ as $(m,3)=1$, and   $m-3>3$ as $m>5$. 
So
$\eta^{\pm 1},\eta^{\pm 2},\eta^{\pm 3}$ are distinct.

\med
(ii) The lemma is true for  if $3|m$,  $(m,5)=1$ and $m\geq 11$.

\med
 Indeed,  $\eta\cdot\eta^4=\eta^5$ and $m-5>5$ for $m\geq11$. So
$\eta^{\pm 1},\eta^{\pm 4},\eta^{\pm 5}$ are distinct. 

\med
(iii) Let $m'$ be the product of all primes dividing $m$. The lemma is true if  $m'<m$.

\med
Indeed, $m'-1,m'+1,2$ are coprime to $m$ so  $\eta^{m'-1},\eta^{-m'-1},\eta^{2}\in\Phi(m)$. As $3m'\leq m$, we have $m\geq 2m'+m'> 2m'+2$.  
If $m'>3$ then $\eta^{\pm(m'+1)}$, $\eta^{\pm(m'-1)}$,
$\eta^{\pm 2}$ are distinct, and $\eta^{m'-1}\eta^{-m'-1}\eta^{2}=1$, as required. Let $m'=3$;
then $m\ge 27$ by assumption. Then $\eta^8\eta^{-10}\eta^2=1$ and $\eta^{\pm 8},\eta^{\pm 10},\eta^{\pm 2}$ are distinct.

\med
(iv) The lemma is true if $m=m'$.

\med
In view of (i), (ii) and (iii) we can assume that  $m>15$, $15|m$ and $(15,m/15)=1$.
 Let $\eta_1=\zeta_3\zeta_5x\in\Phi(m)$, where $\zeta_i\in\Phi(i)$ for $i=3,5$, and $x\in\Phi(m/15)$.
Set $\eta_2=\zeta_3\zeta_5 x^2$ and
$\eta_3=\zeta_3\zeta_5^3x^{-3}$. Then $\eta_2,\eta_3\in\Phi(m)$ and $\eta_{1}\eta_2\eta_{3}=1$
and $\eta_{1},\eta_2,\eta_{3} $ are distinct. 
\enp

\bl{rr1a} Let $m>3$ be odd.
Then either $  \Lambda_3(\Phi(m))=R(m)$ or $m\in\{5,9,15\}$ and  $ \Lambda_3(\Phi(m))=R(m)\setminus 1$.\el

\bp (i) Note that $\Phi(m)\subseteq \Lambda_3(\Phi(m))$. Indeed, if $\zeta\in\Phi(m)$ then there is $\eta\in\Phi(m)$
such that $\zeta\notin\{\eta,\eta\up\}$, so $\eta\cdot\eta\up\zeta=\zeta$. For $m\in\{5,9,15\}$ the lemma follow by inspection.  If $m\notin\{5,9,15\}$
then $1\in \Lambda_3(\Phi(m))$ by  Lemma \ref{333}. This implies the lemma if $m$ is a prime. 

(ii) The lemma is true if $m$ is a prime power. 

Indeed, let $m=p^k>9$, where $p$ is a prime, and $m'=m/p$ if $p>3$ and $m/9$ if $p=3$.  Let $\zeta\in R(m')$.
Suppose first that $p>5$. Then, by Lemma \ref{333}, there are $\eta_1,\eta_2,\eta_3 \in\Phi(m)$ such that $\mu:=\eta_1\eta_2\eta_3\in R(m')$ and $\eta_1,\eta_2,\eta_3 $ are in distinct cosets $rR(m')$, $r\in \Phi (m)$. Set $\nu=\mu\up \zeta$. Then $\eta_1,\eta_2,\eta_3\nu $ are distinct and $\eta_1\eta_2\eta_3\nu =\mu\nu=\zeta$, as required.

Let $p=5$. For $\eta\in\Phi(m)$ set $x=\eta^5\in\Phi(m/5)$ and $\nu=x\up\zeta$. Then  $\eta, \eta^2,\eta^2\nu$ are distinct and $\eta \eta^2\eta^2\nu=x\nu=\zeta$, as required.

Let $p=3$. For $\eta\in\Phi(m)$ set $x=\eta^3\in\Phi(m/3)$.  Since $m>9$,  the elements $\eta,y\eta,y^4\eta$ are distinct for every $y\in \Phi(m/3)$ and their product is $xy^5$. As $\{y^5:y\in \Phi(m/3)\}=\Phi(m/3)$, we conclude that $x\Phi(m/3)\subset \Lambda_3(\Phi(m))$.
By Lemma \ref{b33}, this implies $\Phi(m/3)\cap \Lambda_3(\Phi(m))\neq \emptyset$. 
Moreover, $\Phi(m')\cap \Lambda_3(\Phi(m))\neq \emptyset$ whenever $m'|m$ as $xu\in\Phi(m/3)$
for $u\in R(m/9)$. 
By replacing $\eta$ by $\eta^j$ for integers $j$
coprime to $m$ we conclude that $\Phi(m')\subset \Lambda_3(\Phi(m))$ for every divisor $m'$ of $m$, whence the result.

(iii) Let $m=m_1\cdots m_k$ be the primary decomposition of $m$, and $\zeta=\zeta_1\cdots \zeta_k$, where $\zeta_i\in
R(m_i)$. If $m_i\notin \{3,5,9\}$ for $i=1\ld k$ then the result follows from applying the above to every 
$\zeta_i$,  $i=1\ld k$.

(iv) Suppose that $m_i\in\{3,5,9\}$ for some $i\in\{1\ld k\}$. Then $m\in \{3m',5m',9m',15m',45m'\}$,
where $(m',m/m')=1$. Then $\zeta=\al\be$, where $\beta\in R(m')$  and $\al\in R(m/m')$. 

Suppose first that $m'>1$. Then, by (iii), $\be=\si_1\si_2\si_3$, where $\si_1,\si_2,\si_3\in\Phi(m')$ are distinct. Set $\eta_j=\theta_j\si_j$ for $j=1,2,3$ and $\theta_j\in\Phi(m/m')$. Then $\eta_j\in\Phi(m)$ and $\eta_1,\eta_2,\eta_3$ are distinct 
as so are $\si_1,\si_2,\si_3$. As  $\theta_j\in\Phi(m/m')$ are arbitrary here,    $\al=\theta_1\theta_2\theta_3 $ for some $\theta_1,\theta_2,\theta_3\in\Phi(m/m')$, and hence $\zeta=\eta_1\eta_2\eta_3$ with this choice of $\theta_1,\theta_2,\theta_3$. 

Suppose  that $m'=1$. Then $m=45$ as $m\neq 3,5,9,15$ by assumption.  We can again write $\zeta=\al\be$, where $ \al\in R(5)$, $\be \in R(9)$. As the case with $\zeta\in\Phi(45)\cup 1$ is settled in (i), we have $\al=1, \be\in R(9)\setminus 1$, or $\al\in\Phi(5)$, $\be\in R(3)$. 

Suppose that $\al=1, \be\neq 1$. By (i), $\be=\si_1\si_2\si_3$, where $\si_1,\si_2,\si_3\in\Phi(9)$ are distinct. Set $\eta_1=\om\si_1, \eta_2=\om^2\si_2,\eta_3=\om^2\si_3$ with $\om\in\Phi(5)$. Then  $\eta_1,\eta_2,\eta_3\in\Phi(45)$ are distinct and $\be=\eta_1\eta_2\eta_3$.

Suppose that $\al\neq1 $. Then $\al=\si_1\si_2\si_3$, where $\si_1,\si_2,\si_3\in\Phi(5)$ are distinct. If $\be=1$, set  $\eta_1=\om\si_1, \eta_2=\om^7\si_2,\eta_3=\om\si_3$, where $\om\in\Phi(9)$. 
If $ \be\in\Phi(3) $,  set   $\eta_1=\om\si_1, \eta_2=\om^4\si_2,\eta_3=\om^7\si_3$, where $\om\in\Phi(9)$ and $\om^{12}=\be$. 
Then  $\eta_1,\eta_2,\eta_3\in\Phi(45)$ are distinct and $\al\be=\eta_1\eta_2\eta_3$.\enp

\bl{rr1b} For $m>3$ odd  set $\Lambda^*_3\Phi(m)=\{\eta_1\eta_2\eta_3\up:\eta_1,\eta_2,\eta_3\in\Phi(m)$
 are distinct$\}$. Then  $  \Lambda^*_3\Phi(m)=R(m)$. 
\el

\bp (i) Note that $\Phi(m)\subseteq \Lambda^*_3(\Phi(m))$. Indeed, if $\zeta\in\Phi(m)$ then there is $\eta\in\Phi(m)$
such that $\zeta\notin\{\eta,\eta\up\}$, so $\eta\cdot\eta\up(\zeta\up)\up=\zeta$. For $m\in\{5,9,15\}$ the lemma follow by inspection.  If $m\notin\{5,9,15\}$
then $1\in \Lambda^*_3(\Phi(m))$ by  Lemma \ref{333}. This implies the lemma if $m$ is a prime. 

(ii) The lemma is true if $m$ is a prime power. 

Indeed, let $m=p^k>9$, where $p$ is a prime, and $m'=m/p$ if $p>3$ and $m/9$ if $p=3$. Let $\zeta\in R(m/p)$.
Suppose first that $p>5$. Then, by Lemma \ref{333}, there are $\eta_1,\eta_2,\eta_3 \in\Phi(m)$ such that $\mu:=\eta_1\eta_2\eta_3\up \in R(m')$ and $\eta^{\pm 1}_1,\eta^{\pm 1}_2,\eta^{\pm 1}_3 $ are in distinct cosets $rR(m')$, $r\in \Phi (m)$. 
Set $\nu=\mu \zeta\up $. Then $\eta_1,\eta_2,\eta_3\nu $ are distinct and $\eta_1\eta_2(\eta_3\nu)\up =\mu\nu\up =\zeta$, as required. 

Let $p=5$. For $\eta\in\Phi(m)$ set 
$\nu=\zeta\up$. Then  $\eta, \eta^2,\eta^{3}\nu$ are distinct and $\eta \eta^2(\eta^{3}\nu)\up=\nu\up=\zeta$, as required.

Let $p=3$. For $\eta\in\Phi(m)$ set $x=\eta^3\in\Phi(m/3)$.  Since $m>9$,  the elements $\eta,y\eta,y^2\eta\up$ are distinct for every $1\neq y\in R(m/3)$ and $\eta\cdot \eta y\cdot(y^2\eta\up)\up=xy\up$.
So we conclude that $x(R(m/3)\setminus 1)=R(m/3)\setminus x\subset \Lambda^*_3(\Phi(m))$.
By replacing $\eta$ by $\eta^2$  we conclude that $R(m/3)\subset \Lambda^*_3(\Phi(m))$. By (i), $\Phi(m)\subseteq \Lambda^*_3(\Phi(m))$,  whence the result.

(iii) Let $m=m_1\cdots m_k$ be the primary decomposition of $m$, and $\zeta=\zeta_1\cdots \zeta_k$, where $\zeta_i\in
R(m_i)$. If $m_i\notin \{3,5,9\}$ for $i=1\ld k$ then the result follows from applying the above to every 
$\zeta_i$,  $i=1\ld k$.

 (iv)  Suppose that $m_i\in\{3,5,9\}$ for some $i\in\{1\ld k\}$. Then $m\in \{3m',5m',9m',15m',45m'\}$,
where $(m',m/m')=1$. We have $\zeta=\al\be$, where $\beta\in R(m')$  and $\al\in R(m/m')$. 

Suppose first that $m'>1$. Then, by (iii), $\be=\si_1\si_2\si_3\up$, where $\si_1,\si_2,\si_3\in\Phi(m')$ are distinct. Set $\eta_j=\theta_j\si_j$ for $j=1,2,3$ and $\theta_j\in\Phi(m/m')$. Then $\eta_j\in\Phi(m)$ and $\eta_1,\eta_2,\eta_3$ are distinct 
as so are $\si_1,\si_2,\si_3$. As  $\theta_j\in\Phi(m/m')$ are arbitrary here,    $\al=\theta_1\theta_2\theta_3 \up $ for some $\theta_1,\theta_2,\theta_3\in\Phi(m/m')$, and hence $\zeta=\eta_1\eta_2\eta_3\up$ with this choice of $\theta_1,\theta_2,\theta_3$. 

Suppose  that $m'=1$. Then $m=45$ as $m\neq 3,5,9,15$ by assumption.  We can again write $\zeta=\al\be$, where $ \al\in R(5)$, $\be \in R(9)$. As the case with $\zeta\in\Phi(45)\cup 1$ is settled in (i), we have $\al=1, \be\in R(9)\setminus 1$, or $\al\in\Phi(5)$, $\be\in R(3)$. 

Suppose that $\al=1, \be\neq 1$. By (i), $\be=\si_1\si_2\si_3\up$, where $\si_1,\si_2,\si_3\in\Phi(9)$ are distinct. Set $\eta_1=\om\si_1, \eta_2=\om^2\si_2,\eta_3=\om^3\si_3$ with $\om\in\Phi(5)$. Then  $\eta_1,\eta_2,\eta_3\in\Phi(45)$ are distinct and $\be=\eta_1\eta_2\eta_3\up$.

Suppose that $\al\neq1 $. Then $\al=\si_1\si_2\si_3\up$, where $\si_1,\si_2,\si_3\in\Phi(5)$ are distinct. If $\be=1$, set  $\eta_1=\om\si_1, \eta_2=\om^7\si_2,\eta_3=\om^8\si_3$, where $\om\in\Phi(9)$. 
If $ \be\in\Phi(3) $,  set   $\eta_1=\om\si_1, \eta_2=\om^{-2}\si_2,\eta_3=\si_3$, where $\om\in\Phi(9)$ with $\om^{3}=\be$. 
Then  $\eta_1,\eta_2,\eta_3\in\Phi(45)$ are distinct and $\al\be=\eta_1\eta_2\eta_3\up$.\enp

Recall that $\Lambda_i(\Phi(m))=\{\eta_1\cdots \eta_i:\eta_1\ld \eta_i\in\Phi(m)$ are distinct$\}$.

\bl{222} If $m$ is odd and $\phi(m)>6$ then $\Lambda_3(\Phi(m))\subseteq \Lambda_{i}(\Phi(m))$ for $3<i\leq\phi(m)-3$. \el

\bp Let $\zeta\in \Lambda_3(\Phi(m)) $. Observe that
$\Lambda_i(\Phi(m))=(\Lambda_{\phi(m)-i}(\Phi(m)))\up$. Indeed,  the product of all elements in $\Phi(m)$ equals 1, so if $\eta_1\cdots  \eta_i=\zeta$ for distinct elements
 $\eta_1\ld  \eta_i\in \Phi(m)$ then the product of $\phi(m)-i$ remaining elements of $\Phi(m)$ equals $\zeta\up$. 
So we can assume that  $i\leq \phi(m)/2$. Set $\Phi'(m)=\Phi(m)\setminus \{\eta_1^{\pm 1},\eta_2^{\pm 1},\eta_3^{\pm 1}\}$. Then $\Phi'(m)\up=\Phi'(m)$.
As $\phi(m)-6\geq 2(i-3)$, we can choose  $\eta_{4}\ld \eta_i\in\Phi'(m)$ so that
$\eta_{j+1}=\eta_{j}\up$ for $j\in\{4\ld i-1\}$ even. Then $\eta_1\cdots \eta_i=\zeta$ and 
$\eta_1\ld \eta_i$ are distinct. \enp

\bl{112}  Let $m>3$ odd  and $i\in\{2\ld \phi(m)-2\}$.  $(1)$  $R(m)\setminus 1\subseteq   \Lambda_i(\Phi(m))$.

$(2)$  $ \Lambda_i(\Phi(m))=R(m)$ unless $i$ is odd and $m\in\{5,9,15\}$.
\el

\bp This follows from Lemmas \ref{33a}, \ref{rr1a} and \ref{222}.  \enp

\subsection{The set $\Delta(m)$}

For $m\in \ZZ$, $m>1$ odd, denote by $\mathcal{P}(m)$ the set of all strings $\eta_1\ld \eta_{\phi(m)/2}\in \Phi(m) $ such that $\eta_{1}^{\pm 1}\ld \eta_{\phi(m)/2}^{\pm 1 }$ are distinct. Then 
set $\Delta(m)=\{ \eta_1\cdots \eta_{\phi(m)/2}: \{\eta_1\ld \eta_{\phi(m)/2}\}\in \mathcal{P}(m)\}$. In other words,  $x\in F$ lies in $\Delta(m)$ \ii   $x= \eta_1\cdots \eta_{\phi(m)/2}$ for some $\eta_1\ld \eta_{\phi(m)/2}\in \Phi(m) $ such that $\eta_{1}^{\pm 1}\ld \eta_{\phi(m)/2}^{\pm 1 }$ are distinct.

Note that $x\in \Delta(m)$ implies $x^i\in\Delta(m)$ whenever $(i,m)=1$.

\bl{ni5} Let $m\not\in\{1,3,5,9\}$ be an odd natural number and $k=\phi(m)$. Then $1\in \Delta(m)$. In other words, there are elements $\eta_{1}\ld \eta_{k/2}\in \Phi$ such that $\eta_{1}^{\pm 1}\ld \eta_{k/2}^{\pm 1}$ are distinct and $\eta_{1}\cdots \eta_{k/2}=1$.
\el

\bp (i) Suppose first that $m=p$ is a prime. Then $p>5$.
Set $b=1+2+\cdots +(p-1)/2$ and let $b=lp+t$, where $0\leq t<p$. As $b=\frac{p^2-1}{8}$, we have $t\neq 0$. Suppose first that $t$ is even. Let $b'=b-t$. Then $b'$ can be obtained from $b$ by replacing in the sum $1+2+\cdots +(p-1)/2$ the term $t/2$ by $-t/2. $ Then $b'$ satisfies the conclusion of the lemma. Suppose  that $t$ is odd. Then we replace the term $(p-1)/2$ by $-(p-1)/2 $ to obtain $b+1$ as the new sum. Then we replace the term $(t+1)/2$ by $-(t+1)/2$, which is possible unless $(t+1)/2=(p-1)/2$. If the latter does not hold then   the lemma follows. If   $(t+1)/2=(p-1)/2$ then $t=p-2$ and 
$\frac{p^2-1}{8}=lp+p-2$, whence $p^2-(l+1)p=-15$ and $p\in\{3,5\}$.

(ii) \st $m=p^a>p$.  Let $m=pm'$. Then the mapping $h:R(m)\ra R(p)$ defined by $\zeta\ra \zeta^{m'}$ for $\zeta\in R(m)$ is a group \ho and $h(\Phi(m))=\Phi(p)$.
Suppose first that $p>5$. Then, by (i), there are $\eta_1\ld \eta_{(p-1)/2}\in\Phi(m)$ such that
$h(\eta_1)^{\pm 1}\ld h(\eta_{(p-1)/2})^{\pm 1}$ are distinct, and $h(\eta_1)
\cdots h(\eta_{(p-1)/2})=1$.  

Set $Y_i=\{\eta_ix:x\in R(m')\}$ (so $h(Y_i)=\eta_i^{m'}$) and $Y=Y_1\cup \ldots \cup Y_{(p-1)/2}$. Then $Y\subset\Phi(m)$.
 As $m$ is odd, $y\neq y\up$ for $y\in Y_i$, and obviously $y_i\neq y_j^{\pm 1}$ for
$y_i\in Y_i,y_j\in Y_j$ for $i,j\in\{1\ld (p-1)/2\}$, $i\neq j$. In addition, $\Pi_{y\in Y_i}\,y=\Pi_{x\in R(m')}\,\eta_ix=\eta_i^{m'}$, so $\Pi_{y\in Y}\,y=1$.  This implies the result in this case.

Suppose that   $p\in\{3,5\}$. If $m=27$ then $\phi(m)=18$. Set $\{\eta_1\ld \eta_{9}\}=\{\zeta,\zeta^2,\zeta^{23},\zeta^{5},\zeta^{7},\zeta^{8},\zeta^{10},\\ \zeta^{11},\zeta^{14}\}$ for $\zeta\in\Phi(27)$.
Then $\Pi_{i}\eta_i=\zeta^{81}=1$ and $\eta_i\neq \eta_i\up $,  $\eta_i\neq \eta_j^{\pm 1}$ for  $i,j\in\{1\ld 9\}$, $i\neq j$, whence the lemma in this case. If $27|m$, the result follows as above.

 Suppose  that $m=25$, so $k=20$.  Set $\{\eta_1\ld \eta_{10  }\}=\{\zeta,\zeta^2,\zeta^{3},\zeta^{4},\zeta^{6},\zeta^{18},\zeta^{8},\zeta^{9},\zeta^{11},
\zeta^{13}\}$ for $\zeta\in\Phi(25)$. Then $\Pi_{i}\eta_i=\zeta^{75}=1$ and $\eta_i\neq \eta_i\up $, $\eta_i\neq \eta_j^{\pm 1}$ for  $i,j\in\{1\ld 10\}$, $i\neq j$. If $25|m$, the result follows as above.

(iii) Let $m=p^am'$ with $m'>1$ and $(p,m')=1$. Then $\Phi(m)=\Phi(p^a)\Phi(m')$. Moreover, if $\{y_1\ld y_{\phi(p^a)/2}\}\in\mathcal{P}(p^a)$ and $Y=y_1\Phi(m')\cup\cdots\cup y_{\phi(p^a)/2}\Phi(m')$ then $Y\in\mathcal{P}(m) $ and the product of the elements in $Y$ equals $y_1\cdots y_{\phi(p^a)/2}$. 
If $p^a\neq 3,5,9$ then this product is equal to 1 for some choice of the string $y_1\cdots y_{\phi(p^a)/2}$.   So we are done in this case. 

(iv) We can now  assume $m$ to be coprime to any prime greater than 5, so $m=3^a5^b$ for some $a,b\geq 0$,
and $a\leq 2,b\leq 1$. Therefore, we are left with the cases where $m\in \{15,45\}$, which can be varified by inspection. For instance, if  $m=15$ then $\phi(m)=8$;  set $\{\eta_1\ld \eta_{4 }\}=\{\zeta,\zeta^2,\zeta^{4},\zeta^{8}\}$ for $\zeta\in\Phi(15)$. Then $\Pi_{i}\eta_i=\zeta^{15}=1$ and $\eta_i\neq \eta_j^{\pm 1}$ for $i\neq j$, $i,j\in\{1,2,3, 4\}$ and $\eta_i\neq \eta_i\up $. 
 \enp

\begin{corol}\label{bt4} 
Let $m>1$ be an odd natural number and $k=\phi(m)$. Then $\Phi(m)\subseteq \Delta(m)$.  
\end{corol}

\bp The cases with $m=3,5,9$ follow by inspection. If this does not hold, by Lemma \ref{ni5} $\eta_{1}\cdots \eta_{k/2}=1$ for some $\eta_{1}\ld \eta_{k/2}\in \Phi(m)
$ such that $\eta_{1}^{\pm}\ld \eta_{k/2}^{\pm}$ are distinct, that is,  $\{\eta_{1}, \ld \eta_{k/2}\}$ $\in\mathcal{P}(m)$. In this product replace $\eta_1$ by $\eta_1\up$.
The set $Y:=\{\eta_{1}\up, \eta_2\ld \eta_{k/2}\}$ lies in $\mathcal{P}(m)$ and the product of elements of $Y$ equals $\eta_1^{-2}\in\Phi(m)$. This implies the result as  $\Phi(m)=\{x^i:  (i,m)=1\}$ for every fixed $x\in \Phi(m)$. 
\enp

\bl{n66} Let $m>1$ be an odd natural number and $k=\phi(m)$. Then  $(R(m)\setminus 1)\subseteq\Delta(m)$.
In other words, if $1\neq \zeta\in R(m)$ and   $k=\phi(m)$ then
there are elements $\eta_{1}\ld \eta_{k/2}\in \Phi(m)
$ such that $\eta_{1}\cdots \eta_{k/2}=\zeta$ and $\eta_{1}^{\pm 1}\ld \eta_{k/2}^{\pm 1}$ are distinct.
\el

\bp By  Corollary \ref{bt4}, $\Phi(m) \subset\Delta(m)$. In particular, the lemma is true if $m$ is a prime.

(i) Suppose that $m=p^a, a>1$ for some prime $p$. Set $m'=m/p$. Then $\Phi(m)=\cup_{i=1}^{p-1} \eta_i R(m')$, where $\eta_1\ld \eta_{p-1}$ are representatives of distinct cosets $R(m)/R(m')$.

Suppose first that $p>5$. By Lemma \ref{ni5}, there are $ \eta_1\ld  \eta_{(p-1)/2}\in \Phi(m)$ such that $\eta_1^{\pm m'}\ld\\ \eta_{(p-1)/2}^{\pm m'}$ are distinct and $\eta^{ m'}_1\cdots \eta^{ m'}_{(p-1)/2}=1$. Moreover,  $Y:=\eta_1R(m')\cup\cdots\cup \eta_{(p-1)/2}R(m')\subset \mathcal{P}(m)$ and the product of the elements in $Y$
equals 1. For every $i=1\ld (p-1)/2$ choose an arbitrary element  $\eta_ix_i\in \eta_iR(m')$ and replace it in the product $\Pi_{y\in Y}\,y$ by the inverse of it. Then we obtain the set $Y':=(Y\setminus \{\eta_1x_1\ld \eta_{(p-1)/2}x_{(p-1)/2}\})\cup \{\eta\up_1x\up_1\ld \eta_{(p-1)/2}\up x_{(p-1)/2}\up\}$. The product of the elements of $Y'$ equals
$\eta^{-2}_1\cdots \eta^{-2}_{(p-1)/2}x_1^{-2}\cdots x_{(p-1)/2}^{-2}=(\eta_1\cdots \eta_{(p-1)/2})^{-2}x_1^{-2}\cdots x_{(p-1)/2}^{-2}$. As $\eta_1\cdots \eta_{(p-1)/2}\in R(m')$
and $x_i\in R(m')$ are arbitrary, we conclude that $\Pi_{y\in Y'}\,y$ can be an arbitrary element of $R(m')$
(for a suitable choice of $x_i$'s). As $\Phi(m)\subset\Delta(m)$ by the above, the lemma in this case follows. 

Let $p\leq 5$. One can check that the lemma is true for $m=9,25,27$. 
In addition, for $m\in\{25,27\}$, by Lemma \ref{ni5}, there are $ \eta_1\ld \eta_{\phi(m)/2}\in \Phi(m)$  such that $\eta_1^{\pm 1}\ld \eta_{\phi(m)/2}^{\pm 1}$ are distinct and $\eta_1\cdots \eta_{\phi(m)/2}=1$. Then we  repeat the above reasoning with $m/m'\in\{25,27\}$ in place of $p$
to obtain the result.

(ii) Let $m=p^am'$, where $m'>1$ and $(p,m')=1$. Then $\Phi(m)=\Phi(p^a)\Phi(m')$. Let $l=\phi(p^a)/2=p^{a-1}(p-1)/2$.

Let $\si\in R(p^a)$.  By (i), for every $\si\neq 1$   there are $ \eta_1\ld \eta_{l}\in \Phi(p^a)$ such that 
$\eta_1^{\pm m'}\ld \eta_{l}^{\pm m'}$ are distinct and
 $\eta^{ m'}_1\cdots \eta^{ m'}_{l}=\si$. If $\si=1$  and $p^a>9$ then, by Lemma \ref{ni5},  there are $ \eta_1\ld \eta_{l}\in \Phi(p^a)$ such that $\eta_1^{\pm m'}\ld \eta_{l}^{\pm m'}$ are distinct and  $\eta^{ m'}_1\cdots \eta^{ m'}_{l}=1$.  Then $Y:=\eta_1\Phi(m')\cup\cdots\cup \eta_{l}\Phi(m')\subset \mathcal{P}(m)$. Choose  arbitrary elements  $\eta_ix_i\in \eta_i\Phi (m')$ for $i=1\ld l$ and replace them  in the product $\Pi_{y\in Y}\,y$ by their inverses. Then we obtain the set $Y':=(Y\setminus \{\eta_1x_1\ld \eta_{k}x_{l}\})\cup \{\eta\up_1x\up_1\ld \eta_{l}\up x_{l}\up\}\subset \mathcal{P}(m) $. The product of the elements of $Y'$ equals
$\eta^{-2}_1\cdots \eta^{-2}_{l}x_1^{-2}\cdots x_{l}^{-2}=(\eta_1\cdots \eta_{l})^{-2}x_1^{-2}\cdots x_{l}^{-2}=\si^{-2}x_1^{-2}\cdots x_l^{-2}$. 
 As $x_i\in \Phi(m')$ are arbitrary, by Lemma \ref{t11} we conclude that $\si^{-2}R(m')\subset \Delta(m)$. 
\itf $R(m)=R(p^a)R(m')\subset \Delta(m)$ unless $p^a\in\{3,5,9\} $, $\si=1$ (and then   $(R(p^a)\setminus 1)R(m')\subset \Delta(m)$). 

So the lemma is true if $p|m$ for some prime $p>5$, and if $27|m$ or $25|m$, hence we are left with the cases where $m=3^a5^b$ with $1\leq a\leq 2$ or $b= 1$, that is, with   $m\in\{15,45\}$. In these cases the lemma can be checked straightforwardly. 
  \enp

Then, combining Lemmas \ref{ni5} and \ref{n66}, we have

\begin{corol}\label{c55} $\Delta(m)=R(m)$ if  $m\neq 3,5,9$, and $\Delta(m)=R(m)\setminus \{1\}$ otherwise.   \end{corol}

\subsection{Algebraic groups}

Let $G$ be a simply connected simple algebraic group of rank $n$ defined over an algebraically closed field $F$ of characteristic $p\geq 0$, and let $T$ be a fixed maximal torus of $G$.   Let $\Om={\rm Hom}(T,F^\times)$ be the weight lattice of $G$. We also fix a Borel subgroup of $G$ containing $T$; this defines the fundamental weights  $\om_1\ld \om_n$   of $G$  ordered as in \cite{Bo}. They form a $\ZZ$-basis of $\Om$. The simple roots of $G$ are denoted by $\al_1\ld \al_n$; these are elements of $\Om$. A weight $\om\in\Om$ is called {\it radical} if $\om$
is a $\ZZ$-linear combination of simple roots. The weight 
$\om=\sum a_i\om_i$ ($a_i\in \ZZ$) is called {\it dominant} if $a_i\geq 0$ for all $i$; the dominant weights form a subset $\Om^+$ of $\Om$.  Weight 0 is the one with $a_i=0$ for all $i$. If $\om,\om'\in\Om$ then we write $\om'\preceq \om$ or $\om\succeq \om'$ if the weight $\om-\om'=\sum b_j\al_j$ with $b_j\geq 0$. In \cite{Bo} $\Om$  is embedded into a vector space over the rationals with basis $\ep_1\ld \ep_l$, where $l=8$ if $G$ is of types $E_6,E_7$; $  n+1$ if $G$ is of type $A_n$ and $n$ for the other types; the expressions of the fundamental weights 
$\om_1\ld \om_n$  and simple roots $\al_1\ld \al_n$ in terms of $\ep_1\ld \ep_l$ are tabulated in \cite[Planchees I - III]{Bo}.  

The \ir \reps of $G$ are parameterized by the elements of $\Om^+$, so we can identify them by writing $\rho_\om$ for $\om\in\Om^+$; here $\om$ is called the highest weight of the \rep in question. The group $W=N_G(T)/T$ is called the Weyl group of $G$. The conjugation action of $N_G(T)$ on $T$ induces the action of $W$ on $\Om$. Note that $\rho(T)$ is always completely reducible so $V$, the underlying $G$-module of $\rho$, is a direct sum of $T$-stable subspaces, on each $T$ acts scalarly, and hence define a weight of $G$. The set of such weights is denoted by $\Om(V)$. 
 Note that $\ep_1\ld \ep_l$ can be interpreted as weights of $\rho_{\om_1}$ for groups $G$ of types $A_n,B_n,C_n$ and $D_n$.

 The $G$-module $V$ with highest weight $\om_1$ is called the natural $G$-module for $G$ of classical types; there exists a basis of $V$ under which the elements  $t\in T$ acts on $V$  as diagonal matrices.   This is used to identify elements $g\in G$ in terms of their Jordan normal form. (If $G$ is of types $B_n$ or $D_n$, then $t$
is determined by its matrix on $V$ up to a multiple from $Z(G)$; this is immaterial for our purposes.) As every semisimple element $g\in G$ is conjugate   to an element of $T$,
we can assume that $g\in T$ so $\rho_{\om_1}(g)$ is a diagonal matrix. If $|g|$ is odd then  $g,g'\in T$ are conjugate in $G$ if $\rho_{\om_1}(g), \rho_{\om_1}(g')$ are conjugate in $GL(V)$.

  \med
The \f lemma is often used without an explicit reference. 

\bl{w01} Let G be a semisimple \ag and $g\in G$ a semisimple element. Let V be an \ir G-module. If V has weight $0$ then g has \ei $1$ on V.  
\el

\bp This follows from the standard facts that the maximal tori of $G$ are conjugate and every  semisimple element is contained in a maximal torus. 
\enp

The following lemma is well known; for details see \cite[Lemma 2.12 and Theorem 2.9]{TZ20}:

\begin{lemma}\label{wt1} Let $G$ be a simple \ag of classical type  in defining characteristic $p\geq 0$, and $p\neq 2$ if G is of type $B,C,D$. 
Let $\lambda,\mu\in \Omega^+$  with $\lambda$ $p$-restricted, and $V_\lambda$,
  respectively, $V_\mu$ the associated irreducible $FG$-modules. Then the following facts hold.

 {\rm (1)} If $\mu\prec\lambda$ then    $\Omega(V_\mu)\subseteq \Omega(V_\lambda)$.

{\rm (2)} If $\lambda+\mu$ is $p$-restricted then   $\Omega(V_{\lambda+\mu})=\Omega(V_\lambda\otimes V_\mu)=\Omega(V_\lambda)+\Om( V_\mu)$.

{\rm (3)} If $\lam\neq 0$ is radical then some root is a weight of $V_\lam$; otherwise $\Omega(V_\lam)$ contains some minuscule weight. 

{\rm (4)}  Let $\lambda$ be a $p$-restricted dominant weight. Then
$\Omega(V_\lambda)=\{w(\mu)\ |\ \mu\in \Omega^+, \mu\preceq\lambda,w\in W\}$.
\end{lemma}

Note that minuscule weights are $\om_1\ld \om_{n}$ for type $A_n$, $\om_n$ for type $B_n$ and $\om_1$ for $ C_n$.

 \bl{kn1} \cite[Lemma 5.10(3)]{Z20} Let $ G=SL_n(F)$, $n>1$ odd, and $H=Sp_{n-1}(F)$.
Let V be an \ir $ G$-module. Then $V|_H$ has weight zero.
\el

Recall that an element $g$ of a group $G$ is called {\it real} if $g$ is conjugate to $g\up$.

\begin{corol}\label{kn2} Let $ G=SL_n(F)$, $n>1$ odd, and let $g\in  G$ be  a semisimple element. Let V be a $ FG$-module.  Suppose that  g is real. Then g has \ei $1$ on $V$.
\end{corol}

\bp By \cite[Lemma 5.8]{Z20}, $g$ is contained in a subgroup isomorphic to $ H=Sp_{n-1}(F)$. Then the result follows from Lemma \ref{kn1}.\enp

\bl{la1} For $n\geq 4$ even  let $g=\diag(d_1\ld d_n)\in GL_n(V)$   be a diagonal matrix such that $g$ is odd and  $g\up$ is similar to g.
Suppose that $d_{i_1}d_{i_2}d_{i_3}=1$ for distinct $i_1,i_2,i_3\in\{1\ld n\}$. Then g has \ei $1$ on the $j$-th exterior power   $V_{j}$ of V for $j=2\ld n-2$.
\el

\bp 
The \eis of $g$ on $V_j$ are $d_{i_1}\ld d_{i_j}$ for distinct $i_1\ld i_j\in\{1\ld n\}$.
By reordering of $d_1\ld d_n$ we can assume that $d_{i+1}=d_i\up$ for $i< n$ odd. If $j$ is even then $d_1\cdots  d_j=1$, so we assume $j$ to be odd. Then $n\geq 6$,
and $n=6$ is obvious as well as the cases where $j=3,n-3$ or $1\in\{d_1\ld d_n\}$. Let $n>6, 3<j<n-j$ and $1\notin\{d_1\ld d_n\}$.
Again we can assume that $(i_1,i_2,i_3)=(1,3,5) $. Then $d_1d_3d_5d_7d_8\cdots d_{j+2}d_{j+3}=1$,
whence the result. \enp
  
The \f facts are obvious and we shall use it with no reference.

\bl{3t3} Let G be a  group, and let $g\in G$ be a rational element. Let $\rho$ be a  \rep of $G.$ 

$(1)$ Suppose 
$E( \rho(g))$ contains  a primitive $m$-root of unity for some integer $m$. Then $\Phi(m)\subseteq E (\rho(g))$.

$(2)$  Let $\rho=\rho_1\otimes \rho_2$, where both $\rho_1,\rho_2$ are   \reps of $G.$ Suppose that   $m$ is odd and $\Phi(m)\subseteq E(\rho_i(g))$ for $i=1,2$. Then $R(m)\subseteq E( \rho(g))$. 
\el

\bp (1) Let $\zeta\in E\,( \rho(g)) $ be a primitive $m$-root of unity. As   $g^i$ are conjugate whenever $(i,|g|)=1$, we have $E( \rho(g^i)) = E( \rho(g)) $ and $\zeta^i\in E( \rho(g^i)) $, we are left to observe that $\Phi(m)=\{\zeta^i:(i,|g|)=1\}$. 
  
(2) This follows from Lemma \ref{t11}.
\enp

\bl{gr1} Let G be a simple \ag defined over $\FF_p$ and $g\in G$  a rational \se element. Then g is conjugate to an element of $G(p)$.
\el

\bp Let $Fr$ be the standard Frobenius map on $G$; then, by definition.  $G(p)=G^{Fr}=\{g\in G: Fr(g)=g$. We first observe that the $G$-conjugacy class of $g$ is $Fr$-stable. For this we use a criterion stated in \cite[Ch.II, \S 3, Corollary 3.3]{SpSt}. 

Let $T$ be a maximal torus of $G.$ Then $T$ can be chosen so that $Fr(T)=T$ and $Fr(t)=t^p$ for every $t\in T$.
Let $\rho: G\ra GL(V)$ be an \ir rational \rep of $G$. Then $\rho(T)$ is contained in a maximal torus of  $GL(V)$. Identifying $GL(V)$ with group  of matrices of size $n=\dim V$ with non-zero determinant and replacing $\rho$ by a suitable equivalent representation, we can assume that  $\rho(T)$ consists of diagonal matrices. Then $\rho(Fr(t))=\rho(t^p)=\rho(t)^p$ for every $t\in T$.
 Since $g$ is semisimple, we have $(p,|\rho(g)|)=1=(p,|g|)$.
 Since  $\rho(x)$ is rational in $GL(V)$, it follows that  $\rho(t)$  and 
$\rho(Fr(t))=(\rho(t))^p $ are conjugate in $GL(V)$. As this holds for every $\rho$,  by \cite[Ch.II, \S 3, Corollary 3.3]{SpSt} it follows that $Fr(t)$ and $t$ are conjugate in $G$.
Finally, by Lang's theorem, we conclude that the conjugacy class of $t$ meets $G(p)=G^{Fr}$,
see  \cite[Ch.I, 2.7]{SpSt}. \enp

The \f lemma is sometimes used without an explicit reference. 

\bl{gr8} Let G be a simple \ag 
and let $g\in G$ be  a rational \se element. If $\rho$ is an \irr of G with \hw $\om$.
Then $E(\rho(g))=E(\rho_{p\cdot \om}(g))$. \el

\bp By Lemma \ref{gr1}, we can assume that $g\in G(p)=\{x\in G: Fr(x)=x\}$. In addition, 
 $ \rho_{p\cdot \om}(g)=\rho_{\om}(Fr(g))=\rho_{ \om}(g)$, whence the result. \enp

The \f result is easy and well known.

\bl{gs2} $(1)$ Let $g\in   H=SL_{2n}(F)$ be a real semisimple element. Then g is conjugate to an  element of   $G_1=SO_{2n}(F)$ and $G_2=Sp_{2n}(F)$.

$(2)$ Let $g\in   H=SL_{2n+1}(F)$ be a real semisimple element. Then g is conjugate to an element of $G_1=SO_{2n+1}(F)$.
\el

\bl{zs1} Let $G=Sp_{2n}(F)\subset H=SL_{2n}(F)$. Let $\lam_1\ld \lam_{2n-1}$ be the fundamenttal weights of $H$, and let  $\tau_i$ be an \irr of H with highest weight $\lam_i$, $i\leq n$.

$(1)$ $\rho_{\om_i}$ is a composition factor of $\tau|_G;$ 

$(2)$ if $(i,p)\neq (n,2)$ then every weight of $\tau|_G$ is a weight of $\rho_{\om_i}$.


$(3)$ Let $g\in G$ be a semisimple element. 
Then  $E(\rho_{\om_i}(g)=E(\tau_i(g))$ unless possibly $(i,p)=(n,2)$.
\el

\bp (1) By \cite[Lemma 5.5]{Z20}, every composition factor of $\tau_i|_{G}$ is isomorphic to $\rho_{\om_j}$
for some $j\leq i$ with $i-j$ even. As a Borel subgroup of $G$ is contained in a Borel subgroup of $H$, by  \cite[Lemma 5.5]{Z20} we observe that
$\rho_{\om_i}$ is a composition factor $\tau_i|_{G}$. 

(2) By Lemma \ref{wt1},  $\Om(\rho_{\om_j})\subset \Om(\rho_{\om_i})$,  unless possibly  $p=2$. 
If $i<n$ and $p=2$ then this remains true by \cite[Theorem 15]{z09}. Whence the claim.

(3) follows from (2).\enp

\begin{lemma}\label{abcd} {\rm \cite[Lemma 5.2]{TZ20}} Let $\Omega$ be the root systems of types $A_n$, $n>1$, $B_n$, $n>2$, or $C_n$, $n>2$,  and let $\om\in\Omega$ be a dominant weight. Let $U_1,U_2$ be the subsets of $\Omega$
defined in Table below. 

$(1)$ $\om\succeq \nu$ for some $\nu\in U_1;$

$(2)$ if $\om\notin U_1$ then $\om\succeq \nu$ for some $\nu\in U_2;$

$(3)$ if $\Omega$ is of type $B_n$,  $\om\succ 0$ and $\om\neq \om_1$ then  $\om\succeq \om_2$.
 \end{lemma}

\begin{table}[h]
$$\begin{array}{|l|c|c|}
\hline
\Phi& U_1&U_2 \\
\hline
A_n,n\geq 1&0,\omega_1,\dots,\omega_n& 2\omega_1,2\omega_n,\omega_1+\omega_n,\omega_1+\omega_i,\omega_i+\omega_n, i=2,\dots,n-1\\
\hline
B_n, n\geq2& 0,\omega_n&  \omega_1,\omega_1+\omega_n \\
\hline

C_n, n>2&0,\omega_1&\omega_2,\omega_3\\
\hline
\end{array}$$
\end{table}

\bl{cf1} Let $\Omega$ be a root system  of type $B_n$ and $\om\in\Omega$ a dominant weight. Suppose that $\om\succ \om_n$ and $\om\neq \om_1+\om_n$. Then  $\om\succeq \om_2+\om_n$.
\el

\bp  By Lemma \ref{abcd}, we have $\om\succ \om_1+\om_n$. Let $\om=\sum a_i\om_i$, where $a_i$ are non-negative integers. Note that $0\prec \om_1\prec \cdots \prec \om_{n-1}\prec 2\om_n$ are radical weights. 
As $\om\succ \om_n$, this is not a radical weight and hence $a_n$ is odd. Set $\om'=\om-\om_n$. Then
$\om'\succ 0$. By Lemma \ref{abcd}, either $\om'=\om_1$ (which is excluded by assumption) or $\om'\succeq\om_2$. \itf $\om\succeq\om_2+\om_n$, as claimed. 
\enp

\section{Rational elements in representations of $SL_n(F)$}

In this section $G=SL_n(F),n>1,$ and we identify a semisimple element $g\in G$ with $\rho_{\om_1}(g)$. So we can assume that $g=\diag(d_1\ld d_n)$ is a diagonal matrix of determinant 1.
We view $d_1\ld d_n$ as elements of $F^\times$.
Set  $D=\{d_1\ld d_n\}$.   
  We assume that $|g|$ is odd and
  $g$ is rational, equivalently, $\Phi(|d_i|)\subseteq D$ for every $i=1\ld n$. 


\bl{ts1}
Let $G=SL_n(F),n>4,$ and let $g\in G$ be a semisimple odd order element.
Suppose that $g$ is rational. Then $1$ is not an \ei of $\rho_{\om_3}(g)$ \ii
 one of the \f holds: 

%
$(1)$ $n=6$, $D=\Phi(9)$ and  $|g|=9;$ 

$(2)$ $n=6$, $D=\Phi(5)\cup\Phi(3)$ and $|g|=15;$

$(3)$ $n=8$, $D=\Phi(15)$ and $|g|=15;$ 

$(4)$ $n=8$, $D=\Phi(5)\cup\Phi(3)\cup \Phi(3)$  and $|g|=15;$

$(5)$ $n=10$, $D=\Phi(5)\cup\Phi(9)$ and  $|g|=45;$

$(6)$ $n=14$, $D=\Phi(15)\cup\Phi(9)$  and $|g|=45$.

\noindent In addition, $E(\rho_{\om_1}(g))\subseteq E(\rho_{\om_3}(g))$. 
\end{lemma}
 
\bp We can assume that  $g=\diag(d_1\ld d_n)$ is a diagonal  matrix,  and let $D=\{d_1\ld d_n\}$ be the diagonal entries of $g$. 
The \eis of $\rho_{\om_3}(g)$ are the products $d_{i_1}d_{i_2}d_{i_3}$ for distinct elements  $i_1,i_2,i_3\in \{1\ld n\}$. 

The additional claim is equivalent to saying that  every primitive $|d_i|$-root of unity is an \ei of $\rho_{\om_3}(g)$. This is almost obvious as for every $i$ there are $j,k$ distint from $i$ such that $d_jd_k=1$ and then $d_i=d_id_jd_k\in E(\rho(g))$. In particular, if $d_i=1$ for some $i$ then $1 \in E(\rho_{\om_3}(g))$.  So we can assume that $d_i\neq 1$ for every $i$.

Furthermore, $1\in \Lambda_3(|d_i|)\subset \rho_{\om_3}(g)$ unless $|g_i|\in\{3,5,9,15\}$. 
So the lemma is true unless possibly   $D$ is a union
of subsets $\Phi(3)$, $\Phi(5)$, $\Phi(9)$, $\Phi(15)$, a priory with repetitions. In fact, if a repetition occurs, 
say, $\Phi(m)$ occurs twice then we can try to choose $d_{i_1},d_{i_2}$ in one copy of $\Phi(m)$ and $d_{i_3}$
in another copy of it with $d_{i_3}\up=d_{i_1} d_{i_2}\in\Phi(m)$; this is possible  by Lemma  \ref{33a} 
if $|m|>3$. So $\Phi(m)$ with $m>3$
occurs at most once. In addition,   $\Phi(3)$
occurs at most twice as otherwise we choose  $d_{i_1}=d_{i_2}= d_{i_3}$ in distinct copies of $\Phi(3)$, and then $d_{i_1}d_{i_2} d_{i_3}=1$. 

Next, we rule out the case where $d_j$ is a multiple of $d_i$ and $|d_i|<|d_j|$ for some $i,j\in \{1\ld n\}$. 
Indeed, in this case either $|d_i|=3,|d_j|\in\{9,15\}$ or $|d_i|=5,|d_j|=15$ for somr $d_i,d_j\in D$. Then there are $d_k,d_l\in \Phi(|d_j|)$, $k\neq l$, such that
$d_t:=d_kd_l\in \Phi(|d_i|)$. Then $d_t\up\in \Phi(|d_i|)$ and $d_t\up d_kd_l=1$.

This excludes
all possibilities for $D$ except those listed in the lemma. The fact that $1\notin E(\rho(g))$ for $g$ listed in items $(1)-(6)$     is straightforward.  \enp

\begin{remar}\label{r5} $(1)$ 
As the \reps $\rho_{\om_3}$ and $\rho_{\om_{n-3}}$ are dual to each other,   Lemma {\rm \ref{ts1}} implies a similar result for $i=n-3$.

$(2)$ For arbitrary rational element $g\in G$ primitive $|g|$-roots of unity do not always occur as \eis of $\rho_{\om_3}(g)$. For instance, this happens if $D=\cup_i \Phi(d_i)$ with $i=1,2,3,4$ and $d_1,d_2,d_3,d_4$ are coprime to each other.  
\end{remar}

\bl{de2} Suppose that $g\in SL_n(F)$ is real. Then     $E(\rho_{\om_i}(g))\subseteq E( \rho_{\om_{i+2}}(g))$ for $i\leq j\leq n/2$. \el

\bp Let $\rho_{\om_1}(g)=\diag(d_1\ld d_n)$. Then  $E(\rho_{\om_i}(g))=\{d_{j_1}\cdots d_{j_i}:1\leq j_1< \cdots <  j_i\leq n\}$.
 Let $\zeta\in\Phi(m)$, so $\zeta=d_{j_1}\cdots d_{j_i}$ for some distinct $j_1\ld j_i$.
Let $D'\subseteq D$ be a minimal set such that $d_{j_1}\ld   d_{j_i}\in D'$ and  $D'=(D')\up$.
Then $|D\setminus D'|\geq n-2i$ and $D\setminus D'$ coincides with its inverse. So we can pick 
two elements $x,y\in D\setminus D'$ with $xy=1$. Then $\zeta=d_{j_1}\cdots d_{j_i} =d_{j_1}\cdots d_{j_i}xy\in E( \rho_{\om_{i+2}}(g))$, as required. \enp

\def\hw{highest weight }

\bl{un2} Let $G=SL_n(F)$, n even, let $g\in G$ be a rational odd order element, and  $\rho$  a non-trivial \irr of G.
Then $E(\rho_{\om_1}(g))\subset E(\rho(g))$, unless  n is even,   $\om=p^a\om_j$ with $j$ even, 
and $\rho_{\om_1}(g)$ has exactly two \eis of order multiple to $3$. 
 \el

\bp  Let $\rho_{\om_1}(g)=\diag(d_1\ld d_n)$ and $D=\{d_1\ld d_n\}$. Set $m=|d_i|$ for  $i\in\{1\ld n\}$. We show that $\Phi(m)\subset \rho(g)$ for $m\neq 3$.  Observe that the lemma is true for  $m=1$. 
Indeed, if $n$ is odd then the result is contained in Corollary \ref{kn2}. If $n$ is even then $g$ is contained in a
subgroup $H$ of $G$ isomorphic to $SL_{n-1}(F)$. So the result for $m=1$ follows by applying the result to an \ir constituent of $\rho|_H$. 

\med
From now on we assume that $m>1$. Let $\om$ be the  \hw of $\rho$.  We first consider some special cases. Note that the lemma is trivial for $\om=\om_1$,  and the case with $  \om=\om_{n-1}$ follows by duality reason.

\med
(i) The lemma is true if $\om=\om_j$ for $j$ even. In addition, if $\Phi(3)\subset E(\rho_{\om_1}(g))$ then $xy\in E(\rho(g))$
for all $x\in\Phi(3) $, $y\in\Phi(l)$ whenever $\Phi(l)\subseteq D$.

\med
Suppose first that  $\om=\om_2$. Then the \eis of $\rho(g)$ are $d_id_j$ for $i,j\in\{1\ld n\}$, $i\neq j$.
 If $m>3$ then  the claim follows from Lemma \ref{33a}. If  $m=3$ and  $d_k\in\{1, d_i\}$ for some $k\in\{1\ld n\} $  then the inclusion $\Phi(3)\subset \rho(g)$  is obvious; this holds if $n$ is odd as $d_j=1$ for some $j\in \{1\ld n\}$. 
If there is $m'$ with $3|m'$ and $\Phi(m')\subseteq D\setminus \Phi(3)$ then the claim is again  true  by Lemma \ref{33a}. 
If $\Phi(3)$ occurs once and $\Phi(m')\notin D$ whenever  $3|m'$ then $E(\rho_{\om_2}(g))=\{1\}\cup \{xy\}\cup Y$ for some $Y\subseteq R(|g|/3)$ and all $x\in \Phi(3), y\in D\setminus \Phi(3), $ so $\Phi(3)$ is not contained in $E(\rho_{\om_2}(g))$. This also implies the additional claim. (It is trivial that $1\in E(\rho_{\om_2}(g))$.)

For $j$ even with $2<j\leq n/2$ the result follows by Lemma \ref{de2}. For $\om=\om_{n-j}$ the lemma is true by duality reason.

 \med (iii)  The  lemma  is true if $\om=\om_j$ with $j$ is odd. 

\med
For $j=3$ the result follows from  Lemma \ref{rr1a}, and from Lemma \ref{de2} for $j$ with $3<j\leq n/2$.

\med
(iv) The  lemma   is true if $\om=\om_1+\om_{n-1}$. 

\med
Indeed, the \eis of $\rho(g)$ in this case are $d_id_j\up$ for $i\neq j$, $i,j\in\{1\ld n\}$. So  $\Phi(m)\Phi(m)\subset \rho(g)$ for
$m=|d_i|$. So $R(m)\subseteq E(\rho(g))$ in this case by Lemma \ref{t11}.

\med 
(v) The  lemma is true if $\om$ is $p$-restricted. 

\med
By Lemma \ref{wt1}(3), $\om\succeq\mu$ for some $\mu\in\{\om_1\ld \om_{n-1},0\}$, and if $0\neq\om\succ 0$ then 
$\om\succeq \om_1+\om_{n-1}$. In addition, all weights of $\rho_{\mu}$ are weights of $\rho$, and
if $\mu=0$ then all weights of $\rho_{\om_1+\om_{n-1}}$ are weights of $\rho$. So the result follows from the above,
unless possibly $\om\succ\om_j$ with $j$ even, $\Phi(3)\subset D$ and $\rho_{\om_1}(g)$ has only two \eis of order 3. If $\om\neq \om_j$ here then $\om\succeq \om_1+\om_{j-1}$ or $\om_{j+1}+\om_n$ by Lemma \ref{abcd}. Then  $E(\rho_{\om_1+\om_{j-1}}(g))=
E(\rho_{\om_1}(g))E(\rho_{\om_{j-1}}(g))$ 
and $E(\rho_{\om_{j+1}+\om_n}(g))=E(\rho_{\om_{j+1}}(g))E(\rho_{\om_n}(g))$. Here $j-1,j+1$ are odd,
so, by (iii) above, $E(\rho_{\om_1}(g))\subseteq E(\rho_{\om_{j-1}}(g))\cap  E(\rho_{\om_{j+1}}(g))$. 
Therefore, $R(m)=\Phi(m)\Phi(m)\subseteq E(\rho_{\om_1+\om_{j-1}}(g))\cap E(\rho_{\om_n+\om_{j+1}}(g)) \subset E(\rho(g))$, as requred.

\med
(vi) The lemma  is true if $\om$ is not $p$-restricted. 

\med
By Lemma  \ref{gr8},   we can assume that 
 $\rho_\om$ decomposes as a tensor product of at least two non-trivial factors. Suppose first that $\rho =\rho_{\lam}\otimes \rho_{\mu}$, where $\rho_{\lam}, \rho_{\mu}$ are tensor-indecomposable. Then $\lam=p^a\lam'$, $\mu=p^b\mu'$, where $\lam',\mu'$ are $p $-restricted and $a,b\in \ZZ^+$.
 As the eigenvalues of $V_\lam$ and $V_{\lam'}$ are the same (Lemma \ref{gr8}), we can deal with
$\rho' :=\rho_{\lam'}\otimes \rho_{\mu'}$. If $\Phi(m)\subseteq D$ then, by (v), $\Phi(m)\subseteq E(\rho_{\lam'}(g))$
unless, possibly, $m=3$. If $\Phi(3)\subset E(\rho_{\lam'}(g))$ and $\Phi(3)\subset E(\rho_{\mu'}(g))$ then $\Phi(m)\subset E(\rho_{\lam'}(g))\cdot E(\rho_{\mu'}(g))$. Otherwise, we can assume that $\Phi(3)$ is not contained in $E(\rho_{\lam'}(g))$. Then 
by (i), $xy\in E(\rho_{\lam'}(g))$  for every $x\in\Phi(3)$, $y\in\Phi(l)$ for every $l\neq 3$. By the above, $y\up\in E(\rho_{\mu'}(g))$. Then $x \in
 E(\rho_{\lam'}(g))\cdot E(\rho_{\mu'}(g))$, whence   the result  in this case. If  $\rho_\om$ is a tensor product
of more than two non-trivial factors then we conclude similarly.  \enp

\bp[Proof of Theorem {\rm \ref{tt9}}] 
This is a special case of Lemma \ref{un2}.
\enp

\bl{om12} Let $G=SL_n(F)$, $n\geq 10,$ and let $g=\diag(d_1\ld d_n)\in G$ be an odd order element.
        Let $\rho$ be an \irr of G with \hw $\om_i$, $5\leq i\leq n-5$.
          Suppose that $g$ is rational. Then $1$ is an \ei of $\rho (g)$ unless $n=10$, $D=\Phi(5)\cup\Phi(9)$,     $|g|=45$ and       $i=5$. \el

\bp  Set $m=|g|$.   We can assume that $i\leq n/2$ as the
 \reps $\rho_{\om_i} $ and $\rho_{\om_{n-i}} $ are dual.
 Note that the \eis of $\rho(g)$ are $d_{j_1}\cdots d_{j_i}$, where $j_1\ld j_i\in \{1\ld n\}$ are distinct. 
The case with $i$ even is obvious so we assume $i$ to be odd. In addition, 
by Lemma \ref{de2},  $1\in E(\rho_{\om_3})(g)$ implies 
$1\in E(\rho_{\om_i})(g)$, so we examine the cases in Lemma \ref{ts1}. As $n\geq 10, $ we are left with the cases in items (5),(6) of  Lemma \ref{ts1}.

Suppose that (5) holds.
  We  show that  $1$ is not an \ei of $\rho_{\om_5}(g)$. 
  The \eis of    $\rho_{\om_5}(g)$ are $d_{i_1}d_{i_2}d_{i_3}d_{i_4}d_{i_5}$ for distinct
  $d_1,d_{i_2},d_{i_3},d_{i_4},d_{i_5}\in \Phi(5)\cup \Phi(9)$.
 Suppose that $d_{i_1}d_{i_2}d_{i_3}d_{i_4}d_{i_5}=1$.
Let $s,t$ be the number of terms $d_{i_1}\ld d_{i_5}$ of the left hand side that lie in $\Phi(5)$,$\Phi(9)$, respectively.
By Lemma \ref{rr1a}, $s,t>0$. Obviously, $st\neq 1$. As $|\Phi(9)|=6 $ and $\Pi_{\zeta^k\in\Phi(9)}
 \zeta^k=1$, it follows that $t\neq 5$. Similarly, $s\neq 4$. Then the product of all terms of the product that are  in $ \Phi(9)$
  equals 1 as well as the product of those in $\Phi(5)$.
      Therefore, $(s,t)\in\{(2,3),(3,2)\}$. So we have $d_{i_1}d_{i_2}d_{i_3}=1$
 where the terms are either in $\Phi(5)$ or in $\Phi(9)$.
  By Lemma \ref{rr1a}, this is false.

If (6) holds then let $\zeta\in\Phi(15)$, $\mu\in \Phi(9)$ be such that $\zeta^5=\mu^3$.
Then $\zeta^{-5}=\zeta\up\zeta^{-4}$ and   $\mu.\mu^4.\mu^7=\mu^3 $, so
    $\mu\cdot \mu^4\cdot \mu^7 \cdot \zeta\up \cdot \zeta^{-4}= 1$. Hence the case with $i=5$ is ruled out.
The case with $i=7$ is ruled out by Lemma \ref{de2}. \enp  



  \bp[Proof of Theorem {\rm\ref{th1}}]  The case $\om=0$ is trivial, and will be ignored. 
We can assume that the matrix   $\rho_{\om_1}(g)$ on $V_{\om_1}$ is  diagonal
        and let $D=\{d_1\ld d_n\}$ be the diagonal entries of it regarding the multiplicities. 

 Set $m={\rm max} \{|d_1|\ld |d_n|\}$.
   We can assume that  $\{d_1\ld d_k\}=\Phi(m)$, where $k=\phi(m)$.   In addition, we can reorder $d_1\ld d_n$          and assume that $d_{i+1}=d_i\up$ for $i<n$ odd and $d_n=1$ if $n$ is odd.

      \med
       i) The theorem is true if $n=2$.

\med    In this case $|g|=3$ and $\om=a\om_1$ for some integer $a>1$. The \eis of $\rho(g)$ are
    $\zeta^{a}, \zeta^{a-2}\ld   \zeta^{-a}$, where $\zeta\in\Phi(3)$. As $\zeta^{3}=1$,
    the result easily follows in this case.     So we assume $n>2$ until the end of the proof.

\med
(ii) The theorem is true if $n>3$ and  $\om=\om_i$ for $i\in\{2\ld n-2\}$.

\med
 This follows from Lemmas \ref{ts1} and \ref{om12} for  $i$ is odd and from Lemma \ref{t11} for $i$  even.

          \med 
(iii) The theorem is true if  $\om=\om_2+\om_{n-1}$.

\med
We have $\om_2+\om_{n-1}=2\ep_1+2\ep_2+\ep_3+\cdots +\ep_{n-1}
=\ep_1+\cdots +\ep_n+\ep_1+\ep_2-\ep_n$. 
Then  $(\om_2+\om_{n-1})(g)=(\ep_1+\ep_2-\ep_n)(g)=d_1d_2d_n\up$ as $(\ep_1+\cdots +\ep_n)(g)=1$. 
 By reordering $d_1\ld d_n$ we can assume that $d_1,d_2, d_n\in\Phi(m)$, and, moreover, by
  Lemma \ref{rr1b}, we can assume that $d_1d_2d_n\up=1$ 
(if  $m=3$ then this is straghtforward).

      \med
  (iv) The theorem is true if  $\om$ is $p$-restricted.

\med   

By Lemma \ref{wt1},  $\om\succ \mu$ implies  $\Om(\rho_\mu)\subset\Om(\rho_\om)$, where  $\Om(\rho_\mu),\Om(\rho_\om)$ stand for the set of weights of the respective representations. 
If $\mu=0$ then $\rho$ has weight zero, whence the result in this case. If $\mu\in\{\om_2\ld \om_{n-2}\}$
and $g,n$ are as in the exceptional cases of Lemmas \ref{ts1} and \ref{om12} 
then $\nu(g)=1$ for some weight $\nu$ of $\rho_{\mu}$ by Lemmas \ref{ts1} and \ref{om12}; as $\nu$ is also a weight of $\rho$, the result follows in this cases. 

Suppose that $\mu\in\{\om_1,\om_{n-1}\}$ or $g,n$ are as in one of  the exceptional cases of Lemmas \ref{ts1} and \ref{om12}.
Let  $\si$ be
a    minimal
dominant weight
 such that $\mu\prec\si\preceq\om$. By Lemma \ref{abcd},   if $\mu\in\{\om_1\ld \om_{n-1}\}$  
 then  $\si\in\{\om_1+\om_{i-1},\om_{i+1}+\om_{n-1}\}$ for $1< i<n-1$, 
    and $\si=\om_2+\om_{n-1}$ and $\si=\om_1+\om_{n-2}$ for $i=1,n-1$, respectively. So we have to inspect the following cases: 

\med
(a) $\si=\om_2+\om_{n-1}$ or $\om_1+\om_{n-2}$ (the latter can be ignored   as these two weights are the highest weights of dual modules);

\med
(b) $n>4$, $\si=\om_1+\om_2$ or $\om_4+\om_{n-1}$ 
and one the exceptional cases of Lemmas \ref{ts1} holds;

\med
(c) $\si=\om_1+\om_4$ or $\om_6+\om_{9}$ (the latter can be ignored by the duality reason)  
  and the exceptional case of Lemmas \ref{om12} holds.

\med

Case (a) is considered in (iii), so in this case the theorem is true. 
Consider (b). Suppose first that $\si=\om_1+\om_2 $ and $g$ is as in items $(1)-(6)$ of Lemma \ref{333}. As $\om_1+\om_2=2\ep_1+\ep_2 $,  by reordering of $d_1\ld d_n$ we can assume that
$d_1=\zeta$, $    d_2=\zeta^{-2}$, where $\zeta\in\Phi(5)$ in cases $(1)-(5)$, and $\zeta\in\Phi(15)$ in case $(6)$. Then $(\om_1+\om_{2})(g)=1$, as desired. 

Suppose that $\si=\om_4+\om_{n-1}=2\ep_1+2\ep_2+2\ep_3+2\ep_4+\ep_5+\cdots +\ep_{n-1}=
\ep_1+\cdots+\ep_n+\ep_{1}+\ep_{2}+\ep_{3}+\ep_{4}-\ep_n $, whence $\mu(g)=
(\ep_{1}+\ep_{2}+\ep_{3}+\ep_{4}-\ep_n)(g)$. If $n=6$ then this case can be ignored by the duality reason.
 So we assume that $n\geq 8$  and (5) or (6) of Lemma \ref{ts1} holds. Then $\Phi(9)\subset D$; let  $\zeta\in \Phi(9)$. Then, by reordering of $d_1\ld d_n$, we can assume that
$d_1=\zeta,$ $d_2=\zeta^{2}$, $d_3=\zeta^{4}$, $d_4=\zeta^{7}$, $d_n=\zeta^{5}$. Then
 $\si(g)=\zeta^{1+2+4+7-5}=1$, as desired.

Consider (c).  Here that $n=10$, $g$ is in Lemma \ref{om12} and we can assume that $\si=\om_1+\om_4$.
 We have $\om=\om_1+\om_{4}=2\ep_1+\ep_2+\ep_3+\ep_{4}$. 
 Then we can assume that $d_1=\zeta$,
 $d_2=\zeta^7$, $d_3=\zeta^4$, $d_4=\zeta^5$, where $\zeta\in\Phi(9)$. Then $(\om_1+\om_{4})(g)=(2\ep_1+\ep_2+\ep_3+\ep_{4})(g)=1$.

\med
(v) The theorem is true if  $\om$ is not $p$-restricted.

\med As $1\in \Phi(m)\Phi(m)$, the reasoning in item (vi) of Lemma \ref{un2} shows that the result is true if $|g|>3$.
The case with $|g|=3$ is easy and is left to a reader.\enp 


\def\pri{primitive }
   
\section{The symplectic groups 
}

Throughout this section  $G=Sp_{2n}(F)$ and $H=SL_{2n}(F)$. Let $\om_1\ld \om_n$ be the fundamental weights of $G$ and $\lam_1\ld \lam_{2n-1}$  the fundamental weights of $H$. Both $\om_1\ld \om_n$ and $\lam_1\ld \lam_{2n-1}$ can be viewed as elements of a lattice with basis  $\ep_1\ld \ep_{2n}$ introduced in \cite{Bo} and the expressions of $\om_i$ for $i=1\ld n$ and $\lam_j$ for $j=1\ld 2n-1$ are given in \cite[Planches I, III]{Bo}. In fact, $\ep_1\ld \ep_{2n}$ are functions on the set of diagonal matrices in  $H$ defined by $\ep_i(\diag(t_1\ld t_{2n}))=t_i$ for $i=1\ld 2n$. 
Note that $\ep_1\ld \ep_n$ is a basis of the lattice $\Om(G)$.
We view $G,H$ as the images of the \ir \reps of the respective algebraic groups with \hw $\om_1,\lam_1$, respectively. 

\begin{lemma}\label{rt5} $(1)$ Let $\tau$ be an \irr of $Sp_{2n}(F)$ with \hw $\om_3$,  where $(n,g)$ are
 as in items $(1)-(6)$ of Lemma {\rm \ref{333}}. Then $1$ is not an \ei of $\tau(g)$.

$(2)$ Let $\tau$ be an \irr of $Sp_{10}(F)$ with \hw $\om_5$ and g as
  in item $(5)$ of Lemma {\rm \ref{333}}. Then $1$ is not an \ei of $\tau(g)$.
\el

\bp This follows from Lemma \ref{gs2}.\enp

Let $g\in G$ be a semisimple element. We can assume that $g=\diag(d_1\ld d_{2n})$, where $d_{2n-i+1}=d_i\up$. Then $g$ preserves the non-degenerate bilinear form with Gram matrix $\begin{pmatrix}0&\Id_n\\ -\Id_n&0\end{pmatrix}$, and hence we  can assume that $g\in G\subset H$.  

\bl{un1}  Let $G=Sp_{2n}(F)$,   and let $\rho$ be an \irr of G with \hw $\om_i$, where $1\leq i \leq n$ and $(i,p)\neq (n,2)$. Let $g\in G$ be a rational semisimple element of odd order. Then $E(\rho_{\om_1}(g))\subset E(\rho(g))$ unless $i$ is even
and $g$ is as in the exceptional case of Lemma {\rm \ref{un2}}.
In addition, $1$ is not  an \ei of $\rho(g)$ \ii one of the \f holds:

$(1)$ $i=3$ and $2n,g$ are as in items $(1)-(6)$ of Lemma {\rm \ref{ts1}};

$(2)$ $2n=10$, $i=5$ and $g$ is as in   Lemma {\rm \ref{om12}}.
\el

\bp  
Let  $\tau_i$ be an \irr of $H=SL_{2n}(F)$ with \hw $\lam_i$. By Lemma \ref{zs1},  $E(\rho_{\om_i}(g))\subseteq E(\tau_i(g))$, where the equality holds unless $(i,p)=(n,2)$.  So the result follows in this case from Lemma \ref{un2}. The additional  claim follows from Theorem \ref{th1}.\enp

\bl{od1} Let $p\neq 2$, $G=Sp_{2n}(F)$,    and     let $\rho$ be a \irr of G with \hw $\om\neq 0,p^t\om_1$. Let $g\in G$ be a rational semisimple  element of odd order.

$(1)$ $1\in E(\rho(g))$   unless $\om=p^a\om'$ for some integer $a\geq 0$ and  $(g,2n,\om')$ are as in Lemmas {\rm \ref{ts1}} or {\rm \ref{om12}}.  

$(2)$   $E(\rho_{\om_1}(g))\subset E(\rho(g))$ unless $\om=p^a\om_i$ with $i$ even and g is as in the exceptional case of Lemma {\rm \ref{un2}}. 
\el

\bp If $\om\in\{p^a\om_1\ld p^a\om_{n}\}$ then the lemma follows from   Lemma \ref{un1}.

Suppose first  that $\om$ is $p$-restricted. Then $\om\succ\mu$, where $\mu\in\{\om_1,\om_2\}$ 
and  $E(\rho_\mu)\subseteq E(\rho)$ (Lemma \ref{wt1}). 

(a) Suppose that $\mu=\om_1 $. Then $\Om(\rho_{\om_1})\subset \Om(\rho_{\om})$, whence (2)
in this case. By Lemma \ref{un1},  we can assume that
 $\om\notin\{ \om_1\ld \om_n\}$, and then 
$\om\succeq \om_1+\om_2$ by  \cite[Lemma 2.2, item 5]{Z20}. 
 The weights of $\rho_{\om_1+\om_2}$ are exactly the same as those of $\rho_{\om_1}\otimes\rho_{\om_2}$
(Lemma \ref{wt1}(2)). So $E(\rho_{\om_1+\om_2}(g))=E(\rho_{\om_1}(g))\cdot E(\rho_{\om_2}(g))$.
If $E(\rho_{\om_1}(g))\subset E(\rho_{\om_2}(g))$, then (1) follows as $E(\rho_{\om_1}(g))$
is a union of $\Phi(m)$ for certain divisors $m$ of $|g|$, and $\Phi(m)\Phi(m)=R(m)$
for every $m$ by Lemma \ref{t11}.   Suppose that $E(\rho_{\om_1}(g))$ is not contained in $ E(\rho_{\om_2}(g))$. Then $\Phi(3)\cup \Phi(m)\subseteq E(\rho_{\om_1}(g))$ for some $m$ with
$(3,m)=1$. By Lemma \ref{112}(2),  we have $R(m)\subseteq E(\rho_{\om_2}(g))$, and $1\in \Phi(m)\cdot R(m) \subseteq E(\rho_{\om_1}(g))\cdot  E(\rho_{\om_2}(g))\subset E(\rho(g))   $, as desired. 

(b) Suppose that $\mu=\om_2 $. Then $1\in  E(\rho_{\om_2}(g))$. As  $\om\succ\om_2$ and $\om\notin\{\om_1\ld \om_n\}$, we have $\om\succ 2\om_1$ by \cite[Lemma 2.2(5)]{Z20}.
 The weights of $\rho_{2\om_1}$ are exactly the same as those of $\rho_{\om_1}\otimes\rho_{\om_1}$
(Lemma \ref{wt1}(2)). So $E(\rho_{2\om_1}(g))=E(\rho_{\om_1}(g))\cdot E(\rho_{\om_1}(g))$.
As above, it follows that $E(\rho_{\om_1}(g))\subseteq E(\rho_{2\om_1}(g))\subseteq E(\rho(g)) $, completing the proof for this case and for arbitrary $p$-restricted $\om$. 

 \med
Suppose   that $\om$ is not $p$-restricted. Then $\om=\sum_{i\geq 0} p_i\lam_i$, where $\lam_i$ are $p$-restricted. By Lemma \ref{gr8}, we can assume that the sum here have at least 2 non-zero  terms. As in item (vi) of the proof of Lemma \ref{un2}, we observe that 
$E(\rho(g)) =\otimes_i E(\rho_{\lam_i}(g)) $, and this product consists of at least two terms. 
As above, $ E(\rho_{\lam_i}(g)) $ contains  $E(\rho_{\om_1}(g))$ or  $ E(\rho_{\om_2}(g))$.
We have shown above that  $1\cup  E(\rho_{\om_1}(g))\subseteq  E(\rho_{\om_i}(g))\cdot E(\rho_{\om_j}(g))$ for $(i,j)\in (1,1),(1,2)$.  It is easy to observe that this also true for $(i,j)=(2,2)$. 
\enp

Lemma \ref{od1} implies Theorem \ref{th2} for $p\neq 2$. So we assume until the end of this section that $p=2$.

\bl{o8k} Let    $G=Sp_{2n}(F)$,  $p=2$, and let $g\in G$ be a rational element of odd finite order.
 Let $\rho$ be  an \irr of G of \hw $\om=a_1\om_1+\cdots +a_{n}\om_{n}$.   Suppose that $a_n=0$ and $\om\neq 0,2^l\om_1$ with $l\geq 0$ an integer. Then $1\notin \rho(g)$ \ii cases $(1)$ or $(2)$
  of Lemma {\rm \ref{od1}} holds. 
\el

\bp  By Lemma \ref{gr1}, we can assume that $g\in X=Sp_{2n}(2)$. By \cite[Theorem 1.3]{Z20}, if $g$ is real then
$\rho(g)$ has \ei 1 unless possibly $\lam=2^t\om_i$ for some $i$. For $g$ rational this case is examined in Lemma \ref{un1} for $i<n$ which yields the result in this case.  \enp

Next we consider the case where $p=2$ and $\rho=\rho_{\om_n}$. By Lemma \ref{gr1}, every rational 
element in $G$ is contained in a subgroup $X$ isomorphic to $Sp_{2n}(2)$. Let $V$ be the natural $\FF_2X$-module, and
$V=\oplus_{i=1}^m V_i$, where $V_1\ld V_m$ are minimal non-degenerate $g$-stable subspaces of $V$.
Then  $2d_i:=\dim V_i$ is even, and in some basis of $V$ compatible with the above decomposition
one can write $g=\diag(g_1\ld g_m)$, where $g_i\in Sp_{2d_i}(2)=Sp(V_i)$.
If $g_i$ is \ir in $Sp(V_i)$ then $|g_i|$ divides $2^{d_i}+1$, otherwise $|g_i|$ divides $2^{d_i}-1$.
Following \cite[Theorem 4.5]{Z20}, an index $i$ is called {\it singular} if $|g_i|=2^{d_i}+1$ and  $(|g_i|,|g_j|)=1$ for 
$j\in\{1\ld m\}$ with $j\neq i$, and the number of singular indices is denoted by $Si(g)$.

\bl{si2} Let $g\in X$ be a rational semisimple element. Then $Si(g)\leq 2$. More precisely, if $i$ is a singular index for $g$ then $(d_i,|g_i|)\in\{(1,3),(2,5), (3,9)\}$. In addition,   if   $i\neq j$ are singular indices then 
$5\in\{|g_i|,|g_j|\}$.
\el 

\bp Let $i$ be a singular index for $g$. Then $|g_i|=2^{d_i}+1$ and  $(|g_i|,|g_j|)=1$ for $j\neq i$. This implies   $g_i$ to be rational in $Sp_{2d_i}(2)$. Indeed, if $xgx\up=g^i$ for some $x\in Sp_{2n}(2)$ then $xV_i=V_i$, whence the restriction of $x$ to $V_i$ lies in $Sp_{2d_i}(2)$. (See also \cite[Lemma 3.6]{DoZ}). Therefore, $2d_i\geq \phi(2^{d_i}+1)$. We show that $g_i$  rational  implies $d_i\in\{1,2,3\} $ with $|g_i|=3,5,9$,
 respectively. 

Suppose first that $6\neq 2d_i\geq 4$. By Zsigmondy's theorem \cite[Theorem 5.2.14]{KL}, 
there is a prime $p$ such that $p|(2^{2d_i}-1)$ and coprime to $2^t-1$ for $t<2d_i$. This implies $Sp_{2d_i}(2)$
to have an \ir element $h$ of order $p$, and  $p|(2^{d_i}+1)$. Therefore, $p> 2d_i$. So  $p> 2d_i\geq \phi(2^{d_i}+1)\geq \phi(p)=p-1$. If the latter inequality is strict then $p>  \phi(2^{d_i}+1)\geq 2(p-1)$ as $ \phi(2^{d_i}+1)$ is a multiple of $\phi(p)$, a contradiction.  So $2d_i\geq \phi(2^{d_i}+1)= \phi(p)=p-1$, whence
$2^{d_i}+1=p$ and $2d_i\geq 2^{d_i}=p-1$, which implies $d_i\in\{1,2\}$.   If $d_i=1 $ then $|g_i|=3$ and if $d_i=2$  then $|g_i|=5$.  If $2d_i=6$ then $|g_i|=9$, as claimed. Whence the first statement of the lemma. 

As $(|g_i|,|g_k|)=1$ whenever $i$ is a singular index and $k\neq i$, the additional statement follows. 
\enp

\bl{on9} Let  $p=2$, $G=Sp_{2n}(F)$, and let $g\in G$ be a rational element of odd finite order. 
Let $\rho=\rho_{\om_n}$ be an \irr of G of \hw $\om_n$. 
 Then either $\rho(g)$ has \ei $1$ or one of the \f holds:

$(1)$ $g=\diag(g_1,y)$, where $g_1\in Sp_{2k}(F)$, $k=1,2,3$, $|g_1|=2^k+1$ and
 $y\in Sp_{2n-2k}(F)$ is an arbitrary rational element with $(|y|,|g_1|)=1;$

$(2)$ $g=\diag(g_1,g_2,y)$, where $g_i\in Sp_{2k_i}(F)$, $(k_1,k_2)=(1,2)$ or $(2,3)$,
$|g_i|=2^{k_i}+1$ for $i=1,2$  and $y\in Sp_{2n-2(k_1+k_2)}(F)$ is an arbitrary rational
element with $(|y|,15)=1$.\el

\bp This follows from Lemma \ref{si2} and \cite[Lemma 3.10]{Z20}, where a similar problem is solved for arbitrary semisimple elements in $Sp_{2n}(2)$. To justify, observe that every rational element in $G$ is contained in a
subgroup $X$ isomorphic to $Sp_{2n}(2)$ (Lemma \ref{gr1}) and $\rho|_X$ is irreducible.

If $Si(g)=1$ then we have (1). If $g$ has 2 singular indices $i,j$, say, then $|g_i|,|g_j|$ (up to reordering) are $3,5$ or $5,9$. This is recorded in (2). 
\enp

\medskip {\rm
Let $\om=\sum a_i\om_i\neq\om_n$. In view of Lemma \ref{on9}, we can assume that $a_n\neq 0$. We use    \cite[Proposition  4.9]{Z20}, which is stated   in terms of   $Si(g)$. To be precise, set $\om'=a_1\om_1+\cdots +a_{n-1}\om_{n-1}$;
then we can write $\om'=\sum 2^k\nu_k$, where $\nu_k$'s are 2-restricted weights, and set $\nu=\sum \nu_k$.
Set $\nu=b_1\om_1+\cdots +b_{n-1}\om_{n-1}$ and $\delta(\nu)=\sum b_jj$. Then, by \cite[Proposition  4.9]{Z20}, $\rho_\om(g)$ does  not have \ei 1 \ii   $\delta(\nu)<Si(g)$. Here $\delta(\nu)\geq 1$ as $\nu\neq 0$, so $1\notin E(\rho(g))$ implies   $Si(g)>1$ and also  $Si(g)\leq 2$ by Lemma \ref{si2}. Whence    $Si(g)=2$ and $\delta(\nu)=1$. The letter yields  $\nu=\om_1$ and  $\om'=2^a\om_1$ for some $a\geq 0$, whence  $\om=2^k\om_1+2^t\om_n$. 

Combining this with previous information, we have:

\begin{theo}\label{cn2} Let $G= Sp_{2n}(F)$, $p=2$, and let
 $\rho$ be an \irr of G of \hw $\om=a_1\om_1+\cdots +a_n\om_n$, where $a_n\neq 0$ and $\om\neq p^t\om_n$. Let $g\in G$ be a rational element of odd order. 
Suppose that $1$ is not an \ei of $\rho(g)$. Then   $\om=p^s\om_1+p^t\om_n$, $s,t\geq 0$. Moreover, if V is the natural G-module then  there are g-stable non-degenerate subspaces $V_1,V_2,V_3$ of $V$ such that $V=V_1\oplus V_2\oplus V_3$, the restriction $g_i$ of g on $V_i$ for $i=1,2,3$ are rational elements of $Sp(V_i)$ and one of the \f holds:

$(1)$  $\dim V_1=4$, $|g_1|=5$, $\dim V_2=2$, $|g_2|=3$ and $(|g_3|,15)=1;$

$(2)$  $\dim V_1=4$, $|g_1|=5$, $\dim V_2=6$, $|g_2|=9$ and $(|g_3|,15)=1.$
\end{theo}

\bp The first claim of the theorem has already been established in the reasoning prior to the statement. 
By Lemma \ref{gr1} we can assume that every rational element in $G$ is contained in a subgroup 
$g\in X\cong Sp_{2n}(2)$.  Let $U$ be the natural $\FF_2X$-module. We have also seen above that $Si(g)=2$, so there are exactly two singular indices;  we can assume  these to be $1,2$ and then
  $|g_1|=5$ and $|g_2|\in\{3,9\}$ by Lemma \ref{si2}. Therefore, we can write $g=\diag(g_1,g_2,g_3)$, 
where $(|g_3|,15)=1$.   \itf  there is a $g$-stable orthogonal decomposition $V=V_1\oplus V_2\oplus V_3$, 
such that $\dim U_1=4$ and $\dim U_2\in\{2,6\}$ and  
$g_i$ acts trivially on $V_j$ for $i,j\in\{1,2,3\}, i\neq j$.  So the lemma follows.
\enp

 \bp[Proof of Theorem {\rm \ref{th2}}] The theorem follows from Theorem \ref{od1} for $p\neq 2$,
and from Lemmas \ref{o8k}, \ref{on9} and Theorem \ref{cn2} for $p=2$. 
\enp

\section{The group of type $B_n$}

As $B_2\cong C_2$, we can assume here that $n>2$ and $p\neq 2$. Then universal (simply connected) version of  a simple \ag $G$ of type $B_n$ is known to be $Spin_{2n+1}(F)$, and $\rho_{\om_1}(G)\cong SO_{2n+1}(F)$.

The simple roots of $B_n$ are $\al_i=\ep_i-\ep_{i+1}$ for $i=1\ld n-1$ and $\al_n=\ep_n$. In addition,
$\om_i=\ep_1+...+\ep_i$ for $i<n$ are radical weights and $\om_n=\frac{1}{2}(\ep_1+...+\ep_n)$ is not radical.
Let $T$ be a maximal torus of $G$ and $t\in T$. Then  we can assume that $\rho_{\om_1}(t)=\diag(t_1\ld t_n,1, t_n\up\ld t_1\up)$ with $t_1\ld t_n\in F^\times$, and  $\ep_i(t)=t_i$.

For a diagonal matrix $D$ we denote by $[D]$ the set of diagonal entries of $D$.
Recall that
 $\Delta(m)=\{ x\in F:x=\eta_1\cdots \eta_{\phi(m)/2} \}$, where $\eta_1\ld \eta_{\phi(m)/2}\in \Phi(m) $ are such that $\eta_{1}^{\pm}\ld \eta_{\phi(m)/2}^{\pm}$ are distinct.

 \bl{t66} Let $g\in G=Spin_{2n+1}(F)$ be a rational semisimple element and 
$D=\rho_{\om_1}(g)=\diag(D_1\ld D_r,\Id)\in SO_{2n+1}(F)$, where $D_i$ is a  diagonal $(k_i\times k_i)$-matrix $(i=1\ld r)$ such that $k_i=\phi(m_i)>1$ for some $m_i> 1$ odd, and   $[D_i]= \Phi(m_i) $. 
Then $E(\rho_{\om_n}(g))=\Delta(m_1)\cdots \Delta(m_r)$. 
\el

\bp  We can reorder the diagonal entries of $D$ so that the new matrix $D'$, say, be of shape $\diag(M_1\ld M_r,\Id,M_r\up,\ld, M_1\up)$, where $D_i=\diag(M_i,M_i')$ for $i=1\ld r$. 
Let $d_1\ld d_n$ be the top diagonal entries of $D'$. Then $\ep_j(g)=d_j$ for $j=1\ld n$. The weights of $\rho_{\om_n}$ are $\frac{1}{2}(\pm \ep_1\pm\cdots \pm \ep_n)$. As $|g|$ is odd and $g$ is rational in $G$, the \eis of $\rho(g)$ and of $\rho(g^2)$ are the same, and hence the \eis of $(\frac{1}{2}(\pm\ep_1\pm \cdots \pm\ep_n))(g)$ and of $ (\pm\ep_1\pm \cdots \pm\ep_n)(g)$ are the same. \itf $E(\rho_{\om_n}(g))= (\pm\ep_1\pm \cdots \pm\ep_n)(D)=d_1^{\pm 1}\cdots d_n^{\pm 1}=\Delta(m_1)\cdots \Delta(m_r)$. 
\enp

\bl{pp2} In notation of Lemma {\rm \ref{t66}} we have  
   $E(\rho_{\om_n}(g))=R(|g|)$, unless there is exactly one $i$ with $ m_i\in\{3,5,9\}$, $(m_i,m_j)=1$ for $j\neq i$ and then  $E(\rho_{\om_n}(g))=\Delta(m_i)R(|g|/m_i)$, or there is exactly one $i$ with $ m_i=5$ and exactly one $j$ with $ m_j\in\{3,9\}$, $(m_s,15)=1$ for $s\in\{1\ld r\}$, $s\neq i,j$ and then    $E(\rho_{\om_n}(g))=\Delta(5)\Delta(m_j)R(|g|/5m_j)$.   
\el

\bp By Corollary \ref{c55}, $\Delta(m_i)=R(m_i)$ unless $m_i\in\{3,5,9\}$. One observes that $\Delta(3)\Delta(9)=R(9)$ and $\Delta(m)\Delta(m)=R(m)$ for $m\in\{3,5,9\}$. In addition, $R(m_i)\cdot R(m_j)=R(c)$, and if $(m_i,m_j)>1$ then  $(R(m_i)\setminus 1)\cdot (R(m_j)\setminus 1)=R(c)$, where  $c$ is the least common multiple of $m_i, m_j$ (see Lemma \ref{t2p}). Note that $|g|$ is  the least common multiple of $m_1\ld m_r$.  It follows from Lemma \ref{t66} that $E(\rho_{\om_n}(g))=R(|g|)$, unless there is exactly one $i\in\{1\ld r\}$ with  $m_i\in\{3,9\}$
or exactly one $i$ with $m_i=5$ and exactly one $j\in\{1\ld r\}$  with  $m_j\in\{3,9\}$. In the former case 
$E(\rho_{\om_n}(g))=\Delta(m_i)R(|g|/m_i)$ with $ m_i\in\{3,5,9\}$, in the latter one $E(\rho_{\om_n}(g))=\Delta(5)\Delta(m_j)R(|g|/5m_j)$ with $m_i\in\{3,9\}$.  
\enp

\bl{o10} Let  
$G=Spin_{2n+1}(F)$, let
 $\rho$ be an \irr of G of \hw $\om_n$ and  $g\in G$  an odd order element. Suppose that g is rational. Then  one of the \f holds:

$(1)$ $E(\rho(g))=R(|g|);$
 
 $(2)$ $E(\rho(g))=(R(m)\setminus 1) R(|g|/m)=(R(|g|)\setminus R(|g|/m)$ with $m\in\{3,5,9\}$ and $(m,|g|/m)=1;$

$(3)$  
$E(\rho(g))=(R(m)\setminus 1)\Phi(5) R(|g|/5m)$ with $m\in\{3,9\}$ and $(5m,|g|/5m)=1;$
in other words, $E(\rho(g))$ consists of all elements of $R(|g|)$ whose orders are multiples of $15.$

\med
In addition, in notation of Lemma {\rm \ref{t66}} let 
 $D=\diag(D_1\ld D_r,\Id)=\rho_{\om_1}(g)$ and let $m_i$ be as in Lemma $\ref{pp2}$.  Then
$(2)$ holds \ii there is exactly one $m_i\in\{3,5,9\}$ with $[D_i]=\Phi(m_i)$  and $(m_i,m_j)=1$ for $j\neq i$ $(i,j\in\{1\ld r\}  )$, and  $(3)$ holds \ii there is exactly one $m_i\in\{3,9\}$ with  $[D_i]=\Phi(m_i)$, exactly one $j$ with  $[D_j]=\Phi(5)$, $i,j\in\{1\ld r\}$, and $(5m_i,m_k)=1$ for $k\neq i , j$.
\el

\bp  
By Corollary \ref{c55}, $\Delta(m_i)=R(m_i)$ unless $m_i\in\{3,5,9\}$ and $\Delta(m)=(R(m)\setminus 1)$ for $m\in\{3,5,9\}$. One observes that $\Delta(3)\Delta(9)=R(9)$ and $\Delta(m)\Delta(m)=R(m)$ for $m\in\{3,5,9\}$. In addition, $R(m_i)\cdot R(m_j)=R(c)$, and if $(m_i,m_j)>1$ then  $(R(m_i)\setminus 1)\cdot (R(m_j)\setminus 1)=R(c)$, where  $c$ is the least common multiple of $m_i, m_j$ (see Lemma \ref{t2p}). Note that $|g|$ is  the least common multiple of $m_1\ld m_r$.  It follows from Lemma \ref{t66} that $E(\rho_{\om_n}(g))=R(|g|)$, unless there is exactly one $i\in\{1\ld r\}$ with  $m_i\in\{3,5,9\}$
or exactly one $i$ with $m_i=5$ and exactly one $j\in\{1\ld r\}$  with  $m_j\in\{3,9\}$. In the former case 
$E(\rho_{\om_n}(g))=\Delta(m_i)R(|g|/m_i)$ with $ m_i\in\{3,5,9\}$, in the latter one $E(\rho_{\om_n}(g))=\Delta(5)\Delta(m_j)R(|g|/5m_j)$ with $m_i\in\{3,9\}$.  


For more details, let $V$ be the underlying  space  for $\rho_{\om_1}$; then 
for every  semisimple element of $G$ the matrix of $\rho_{\om_1}(g)$ 
is conjugate to $\diag(d_1\ld d_n,1,d_n\up\ld d_1\up)$  under a suitable  basis of $V$. The set of diagonal matrices (under this basis) forms a maximal torus of $G$ and  $\ep_i(g)=d_i$ for $i=1\ld n $. As $g$ is rational, $E(\rho(g))=E(\rho(g^2))$ for arbitrary  \rep $\rho$ of $G$.

By Lemma \ref{t66},  the values of $(\frac{1}{2}(\pm\ep_1\pm \cdots \pm\ep_n))(g^2)$ and of $ (\pm\ep_1\pm \cdots \pm\ep_n)(g)$ are the same.   So we determine  
 $ (\pm\ep_1\pm \cdots \pm\ep_n)(g)$, which is  are $d_1^{\pm 1}\cdots d_n^{\pm 1}$, where $d_i\neq d_j^{\pm 1}$ for $i\neq j$, $i,j\in\{1\ld n\}$. 
\enp

\begin{theo}\label{or1}  Let   
$G=Spin_{2n+1}(F)$, $p\neq 2$, 
 let $\rho$ be an \irr of G of \hw $\om$ and $g\in G$ an element of finite order. Suppose that g is rational and $(2,|g|)=1$. Then  one of the \f holds:

$(1)$ $\rho(g)$ has \ei $1;$

$(2)$ $\om=p^j\om_n$ $(j\geq 0)$, and $(2)$ or $(3)$ of Lemma $\ref{o10}$ holds;

$(3)$ $\om=p^i\om_1+p^j\om_n$ $(i,j\geq 0)$, and $(3)$ of Lemma $\ref{o10}$ holds.
\end{theo}

\bp Suppose that (1) does not hold. If $\om=p^j\om_n$ then the result follows from   Lemma \ref{o10}. Suppose that 
(1) holds and  $\om\neq p^j\om_n$.

$(i)$ Suppose that $\om$ is $p$-restricted. Then either $\om\succ 0$ or $\om\succ  \om_n $. In the former case, by Lemma \ref{wt1}, $\rho$ has weight zero, and the result follows by Lemma \ref{w01}. 

Suppose that $\om\succ \om_n$. Then, by Lemma \ref{o10}, $\om_n\prec\lam\preceq \om$, where $\lam=\om_1+\om_n$. 

\med
(a) Suppose first that $\om=\lam$. Note that $E(\rho(g))=E(\rho_{\om_1}(g))\cdot E(\rho_{\om_n}(g))$ by Lemma \ref{wt1}. We use Lemma \ref{o10}. If $(1)$ of  Lemma \ref{o10} holds then $E(\rho(g))=R(|g|)$ and the claim is obvious.  Suppose that (2) of  Lemma \ref{o10} holds. We  show that $E(\rho(g))=R(|g|)$.  
Indeed, we have $\Phi(m)\subset E(\rho_{\om_1}(g))$
 and $E(\rho_{\om_n}(g))$ coincides with
 $\Phi(m)\cdot R(|g|/m)=R(|g|)\setminus R(|g|/m)$, where $m\in \{3,5,9\}$ and $(m,|g|/m)=1$.  As $\Phi(m)\cdot\Phi(m)=R(m)$, it follows that  $R(m)\cdot R(|g|/m)\subset E(\rho(g))$. As $(m,|g|)=1$, we have
$R(m)\cdot R(|g|/m)=R(|g|)$. 

 Suppose that (3) of  Lemma \ref{o10} holds. Then $E(\rho_{\om_n}(g))=\Phi(5)\cdot (R(m)\setminus 1)\cdot R(|g|/5m)$, where
 $m\in\{3,9\}$ and $(15,|g|/5m)=1$. In other words, these are all elements of $R(|g|)$ whose order are multiples of $15$. The spectrum of $\rho_{\om_1}(g)$ coincides with $\Phi(5)\cup \Phi(m)\cup  X$, where $X\subseteq R(|g|/5m)$.
We have $$(\Phi(5)\cup \Phi(m))\cdot \Phi(5)\cdot (R(m)\setminus 1)=(R(5)\cup \Phi(5m))\cdot (R(m)\setminus 1))=$$ $$=(R(5m)\setminus R(5))\cup (\Phi(5m)\cdot (R(m)\setminus 1))=(R(5m)\setminus R(5))\cup (\Phi(5)\cdot R(m)\setminus \Phi(5m))= R(5m)\setminus 1$$ since  the elements of order 5 occurs  in $\Phi(5)\cdot R(m)\setminus \Phi(5m)$. So    $E(\rho(g))=R(|g|)\setminus R(|g|/5m)$.
In particular, $1\notin E(\rho(g))$ and  $\Phi(5)\cdot R(|g|/5m)\subset  E(\rho(g))$.

\med
(b) Suppose that $\om\succ \lam$. Then $E(\rho(g))=|g|$. 

Indeed, by Lemma \ref{cf1}, $\om\succeq \om_2+\om_n$. Then $E(\rho_{\om_2}(g))\cdot E(\rho_{\om_n}(g))\subseteq E(\rho(g))$ by Lemma \ref{wt1}.  If $E(\rho_{\om_n}(g))=R(|g|)$ then $E(\rho(g))=R(|g|)$, and we have (1). Suppose that  $E(\rho_{\om_n}(g))\neq R(|g|)$; then $(2)$ or $(3)$  of  Lemma \ref{o10} holds.

As $n>2$, the weights of $\rho_{\om_2}$ are $\pm\ep_1\pm\ep_2\cup 0$.
Suppose that (2) of  Lemma \ref{o10} holds. By Lemma \ref{pp2}, $E(\rho_{\om_n}(g))=\Phi(m)R(|g|/m)$,  $m\in\{3,5,9\}$, and then
  $E(\rho_{\om_2}(g))$ contains the set $ \Phi(m)X$, where $X$ is the set of the \eis of $\rho_{\om_1}(g)$ whose orders are comprime to $m$. As $X\subset R(|g|/m)$, by Lemma \ref{t11} we have $E(\rho_{\om_2}(g))\cdot E(\rho_{\om_n}(g))=
R(m)\cdot R(|g|/m)=R(|g|)$ as here $(m,|g|/m)=1$.

Suppose that (3) of  Lemma \ref{o10} holds.   By Lemma \ref{pp2}, 
$E(\rho_{\om_n}(g))=(R(m)\setminus 1)\Phi(5)R(|g|/m)$,  $m\in\{3,9\}$.   
As $\Phi(5m)=\Phi(5)\Phi(m) \subset E(\rho_{\om_2}(g))$, we have
 $$R(5m) =\Phi(5m)\Phi(5m)\subseteq   \Phi(5m)(R(m)\setminus 1)\Phi(5)= \Phi(5)\Phi(m)(R(m)\setminus 1)\Phi(5)$$ by Lemma \ref{t11},
and then 
$E(\rho_{\om_2}(g))\cdot E(\rho_{\om_n}(g))$ contains 
$ R(|g|)$ as desired.

\med

$(ii)$ Suppose that $\om$ is not $p$-restricted, say, $\om=\sum_{i\geq 0} p^i\mu_i$, where $\mu_i$'s are $p$-restricted. 
In view of Lemma  \ref{gr8} we have $E(\rho(g))=E(\otimes_i \rho_i(g))$, where each $\rho_i$ is a $p$-restricted \irr of $G$ with \hw $\mu_i$. If  $1\in E(\rho_i(g))$ for  every $i$, then $1\in  E(\rho(g))$. Suppose that $1\notin  E(\rho(g_i))$ for some $i$; we can assume that $i=0$ here.  
For our purpose we can assume that $\mu_i\neq 0$ for every $i$. By Lemma \ref{abcd},
either $\mu_i\succeq\om_n$ or $\mu_i\succeq\om_1$; by the above we can assume that $\mu_0\succeq\om_n$. 

Suppose first that $\mu_j\succeq\om_n$ for some $j>0$, and we can assume that $j=1$.
Then we show that  $E(\rho_0(g))E(\rho_1(g))=R(|g|)$. 
By Lemma \ref{wt1}, $E(\rho_{\om_n}(g))\subseteq E(\rho_0(g))\cap E(\rho_1(g))$, so it suffices to observe that 
$E(\rho_{\om_n}(g))E(\rho_{\om_n}(g))=R(m)$. By Lemma \ref{o10}, we have to inspect the cases (2),(3) of 
that lemma. Observe that $g$ is fixed here, so one can immediately conclude from the expressions for $E(\rho_{\om_n}(g))$ in Lemma \ref{o10} that the required equality holds. 

So we are left with the case where $\mu_j\succeq \om_1$ for all $j>0$.  If $\mu_i\neq \om_1$ 
then, by (i)(b) above,  $E(\rho_{\mu_1}(g))=R(|g|)$, whence $E(\rho(g))=|g|$.

Suppose that $\mu_i=\om_1$ for $i>0$.  We show that $i\leq 1$. 

Indeed, suppose that $\mu_2=\om_1. $ As $p>2$ here, we have $E(\rho_{\om_1}(g))\cdot E(\rho_{\om_1}(g))=E(\rho_{2\om_1}(g))$ and $E(\rho_{\om_n}(g))\cdot E(\rho_{2\om_1}(g))=E(\rho_{2\om_1+\om_n}(g))$ by Lemma \ref{wt1}. Then  $1\in    E(\rho_{2\om_1+\om_n}(g))=R(|g|)$ by (i)(b), whence the result.
\enp

 \bp[Proof of Theorem {\rm \ref{th3}}] This is a self-contained form of Theorem \ref{or1}.
\enp

{\bf Acknowledgement} 
This work has been completed during the author's visit to the Isaac Newton Institute for Mathematical Sciences, Cambridge (June 26th - July 29th, 2022), for his participation in the research programme `Groups, representations and applications: new perspectives'.   I would like to thank the Isaac Newton Institute for  hospitality during my visit and  acknowledge  financial support,   EPSRC grant no. EP/R014604/1, very   encouraging scientific activity at the Institute including seminars and lectures on the latest progress in the area of the programme.


\begin{thebibliography}{1245}

\bibitem{Bo} N. Bourbaki, {\it Groupes et algebres de Lie}, ch. IV-VI, Masson, Paris, 1981.

 
\bibitem{Cu}  J. Cullivan, Fixed-point subgroups of $GL_3(q)$, {\it J. Group Theory} 22(2019), $893 - 914$.

\bibitem{CZ} J. Cullinan  and A.E. Zalesski,  Unisingular Representations in Arithmetic and Lie Theory, 
{\it  European J. Math}. 7(2021), no.4, $1645 - 1667$.
 

\bibitem{DoZ} M. Dokuchaev and  A.E. Zalesski, On the automorphism group of rational group algebras of finite groups, {\it Contemp. Math.} 668(2017), $33 - 51$.

\bibitem{EZ} L. Emmett and A.E. Zalesski, On regular orbits of elements of classical groups in their
permutation representations, {\it Comm. Algebra} 39(2011), $3356 - 3409$.

\bibitem{GT} R. Guralnick and Pham Huu Tiep, Finite simple unisingular groups of Lie type, {\it J. Group Theory} 6(2003), 271 $-$ 310.


\bibitem{HZ09} G. Hiss and A.E. Zalesski, The Weil-Steinberg character of finite classical groups
    with an appendix by Olivier Brunat, {\it Represent. Theory} 13(2009), 427 - 459. Corrigendum: {\it Represent. Theory} 15(2011), $729 - 732$.


\bibitem{Hum} J. Humphreys,  {\it Modular representations of finite groups of Lie type}, Cambridge Univ. Press, Cambridge, 2005.

\bibitem{Hu} B. Huppert, Singer-Zyklen in klassicshen Gruppen,   {\it Math. Z.} 117(1970), 141 $-$ 150.


\bibitem{KL} P. Kleidman and M. Liebeck,  {\it Subgroup structure of classical groups}, Cambridge Univ. Press, Cambridge, 1990.  
 

\bibitem{SpSt} T. Springer and R. Steinberg, Conjugacy classes, In: {\it Seminar on algebraic groups and related finite groups}, Lect. Notes in Math. vol. 131, Springer-Verlag, Berlin, 1970, pp. $167 - 266$.

 



\bibitem{TZ20} D. Testerman and A.E. Zalesski, Spectra of non-regular elements in irreducible representations of simple algebraic groups, {\it North-Western European J. Math}. 7(2021), 7(2021), $185-212.$
 
\bibitem{Wi} R. Wilson, Certain representations of Chevalley groups over $GF(2^n)$, {\it Comm. Algebra} 3(1975), 319 -- 364.

\bibitem{Z88}  A.E. Zalesski, Eigenvalues of matrices of complex representations of finite groups of Lie type. In: {\it Algebra, some current trends}, Lecture Notes in Math. (Springer-Verlag), vol. 1352, 1988'', $206 - 218$.

\bibitem{Z91}  A.E. Zalesski, The \ei $1$ of matrices of complex \reps of finite Chevalley groups, {\it Proc. Steklov Inst. Math}   (1991), no.4, $109 - 119$.

\bibitem{z09}  A.E. Zalesski,  On eigenvalues of group elements in \reps of simple
algebraic groups and finite Chevalley groups. {\it Acta Appl.
Math.} 108(2009), 175 -- 195.
 


\bibitem{Z17}  A.E. Zalesski, Singer cycles in 2-modular \reps of $GL_n(2)$,  {\it Archiv der Math.}
110(2018), $433 - 446$.

\bibitem{Z20} A.E. Zalesski,   
Unisingular representations of finite symplectic groups, {\it Comm. Algebra} 50(2022), $1697 - 1719$. 
 
\end{thebibliography}
\end{document}